\documentclass{article}

\usepackage{amsmath,graphicx,dsfont}

\def\R{\mathbb R}

\def\N{\mathbb N}

\def\L{\mathcal{L}}

\newcommand\Co{\mathbb{C}}
\newcommand\C[1]{\mathcal{C}^{#1}}

\def\tu{\tilde{\upsilon}}

\def\su{\upsilon^*}

\DeclareMathOperator{\im}{Im}

\def\epsilon{\varepsilon}

\def\lp {\left( }
\def\rp {\right) }

\def\ds{\displaystyle}

\newcommand{\be}{\begin{equation}}
\newcommand{\ee}{\end{equation}}
\newcommand{\baa}{\begin{array}}
\newcommand{\eaa}{\end{array}}
\newcommand{\ba}{\begin{eqnarray}}
\newcommand{\ea}{\end{eqnarray}}

\newtheorem{defi}{Definition}
\newtheorem{theo}{Theorem}

\newtheorem{lem}{Lemma}

\newtheorem{prop}{Proposition}

\newtheorem{remark}{Remark}

\newcommand{\Proof}[1]{\noindent\textit{Proof#1. }}
\newcommand{\carre}{\hfill$\Box$\par\addvspace{4mm}}

\usepackage{amsmath,graphicx,dsfont}
\usepackage{epsfig}
\usepackage{psfrag}
\usepackage[english]{babel}
\usepackage{amsmath}
\usepackage{amssymb,subfigure}
\usepackage{mathrsfs}
\usepackage{color}
\usepackage{textcomp}
\usepackage{pstricks,pst-text,pst-3d}
\usepackage{pst-node}
\usepackage{pst-plot}
\usepackage{lmodern}
\usepackage{fontenc}

\begin{document}

\title{Inside dynamics of pulled and pushed fronts}

\author{Jimmy Garnier$^{\hbox{ \small{a,b}}}$, Thomas Giletti$^{\hbox{ \small{b}}}$, Fran\c{c}ois Hamel$^{\hbox{ \small{b,c}}}$ \& Lionel Roques$^{\hbox{ \small{a}}}$\thanks{The authors are supported by the French ``Agence Nationale de la Recherche" within the projects ColonSGS and PREFERED.} \\
\footnotesize{$^{\hbox{a }}$UR 546 Biostatistique et Processus Spatiaux, INRA, F-84000 Avignon, France}\\
\footnotesize{$^{\hbox{b }}$Aix-Marseille Universit\'e, LATP (UMR 6632), Facult\'e des Sciences et Techniques}\\
\footnotesize{Avenue Escadrille Normandie-Niemen, F-13397 Marseille Cedex 20, France}\\
\footnotesize{$^{\hbox{c }}$ Institut Universitaire de France}
}

\date{}
\maketitle

\begin{abstract}\noindent{}We investigate the inside structure of one-dimensional reaction-diffusion traveling fronts. The reaction terms are of the monostable, bistable or ignition types. Assuming that the fronts are made of several components with identical diffusion and growth rates, we analyze the spreading properties of each component. In the monostable case, the fronts are classified as pulled or pushed ones, depending on the propagation speed. We prove that any localized component of a pulled front converges locally to $0$ at large times in the moving frame of the front, while any component of a pushed front converges to a well determined positive proportion of the front in the moving frame. These results give a new and more complete interpretation of the pulled/pushed terminology which extends the previous definitions to the case of general transition waves. In particular, in the bistable and ignition cases, the fronts are proved to be pushed as they share the same inside structure as the pushed monostable critical fronts. Uniform convergence results and precise estimates of the left and right spreading speeds of the components of pulled and pushed fronts are also established.
\end{abstract}


\section{Introduction}

In this paper, we explore the spatial structure of traveling wave solutions of some reaction-diffusion equations. Namely, we consider the following one-dimensional reaction-diffusion model:
\begin{equation}\label{eq:RD}
    \partial_t u(t,x) = \partial^2_{x}u(t,x) +  f(u(t,x)),  \ \ t>0,\ x\in\R,
\end{equation}
where $u(t,x)\in[0,1].$
This equation arises in various scientific domains of application, namely, population dynamics where the unknown quantity $u$ typically stands for a population density~\cite{Bri86,Cancos03,Gri96,Mur02,ShiKaw97}, chemistry~\cite{BilNee91,Fif79b}, and combustion~\cite{BerLar91}. In the context of population dynamics, $f(u)$ corresponds to the population's growth rate.
The nonlinear growth term $f$ in \eqref{eq:RD} is assumed to satisfy
\begin{equation}\label{hyp:f}
f \in C^1([0,1]), \hspace{0.3cm}  f(0)=f(1)=0, \hspace{0.3cm}  \displaystyle  \int_{0}^{1}\!\!{f(s)\,ds}>0
\end{equation}
and to be either of the monostable, bistable or ignition type:
\begin{description}
\item[(A) Monostable] $f$ is \emph{monostable} if it satisfies~\eqref{hyp:f},  $f'(0)>0,$ $f'(1)<0$ and $f>0$ in $(0,1).$\par
In this case, the growth rate $f(u)$ is always positive on $(0,1)$. A classical monostable example is $ f(u)=u(1-u)(1+au),$ with $a\geq0.$ If $a\leq1,$ the per capita growth rate $f(u)/u$ (defined by $f'(0)$ at $u=0$) is decreasing over the interval~$[0,1]$ which means that $f$ is of the KPP (for Kolmogorov-Petrovsky-Piskunov~\cite{KPP}) type. On the other hand, if $a>1,$ the maximum of the per  capita growth rate is not reached at $u=0.$ In population dynamics, this corresponds to a so-called weak Allee effect~\cite{Tur98}.

\item[(B) Bistable]  $f$ is \emph{bistable} if it satisfies~\eqref{hyp:f},  $f'(0)<0,$ $f'(1)<0$ and there exists $\rho\in(0,1)$ such that  $f<0$ in $(0,\rho)$ and $f>0 $ in $(\rho,1).$\par
This hypothesis means that the growth rate $f(u)$ is negative at low densities, which corresponds to a strong Allee effect~\cite{LewKar93,Tur98}. The parameter $\rho$ corresponds to the so-called ``Allee threshold" below which the growth rate becomes negative. For instance, the cubical function $f(u)=u(1-u)(u-\rho),$ where $\rho$ belongs to $(0,1/2)$, is a bistable nonlinearity.

\item[(C) Ignition]  $f$ is of \emph{ignition} type if it satisfies~\eqref{hyp:f}, $f'(1)<0$ and there exists $\rho\in(0,1)$ such that  $f=0$ in $(0,\rho)$ and $f>0 $ in $(\rho,1).$ This reaction term occurs in combustion problems, where $u$ corresponds to a temperature and $\rho$ is the ignition temperature \cite{BerLar91,BerNic85,Roq05}.
\end{description}

The equation \eqref{eq:RD} has been extensively used to model spatial propagation of elements in interaction, in parts because it admits traveling wave solutions. These particular solutions keep a constant profile $U$ and move at a constant speed $c.$ Aronson and Weinberger~\cite{AroWei75,AroWei78} and Kanel'~\cite{Kan61} have proved that equation~\eqref{eq:RD} with monostable, bistable or ignition nonlinearities admits traveling waves solutions of the form $u(t,x)=U(x-ct),$ where $c\in\R$ and the profile $U$ is a $\C3(\R)$ function which satisfies the following nonlinear elliptic equation:
\begin{equation}\label{eq:pb_front}
\left\{\begin{array}{l}
       U''(y)+c\,U'(y)+f(U(y))=0,\ \ y\in\R ,\vspace{3pt} \\
       U(-\infty)=1, \  \ U(+\infty)=0 \ \hbox{ and } \ 0<U<1 \hbox{ on } \R.
      \end{array}\right.
\end{equation}
On the one hand, if $f$ is of the monostable type~$(A),$ there exists a minimal speed $c^*\geq 2\sqrt{f'(0)}>0$ such that equation~\eqref{eq:pb_front} admits a solution if and only if $c \geq c^*$. The solution associated to the minimal speed $c^*$ is called the \emph{critical front},  while those associated to speeds $c>c^*$ are called \emph{super-critical fronts.} On the other hand, if $f$ is of the bistable~$(B)$ or ignition~$(C)$ types there exists a unique speed $c>0$ such that equation~\eqref{eq:pb_front} admits a solution. In all cases,  if the front~$(c,U)$ exists, the profile~$U$ is a decreasing function and it is unique up to shift (see e.g. \cite{AroWei78,BerNic85,FifMcL77}). The asymptotic behavior of $U(y)$ as $|y|\to+\infty$ is also known~(see Section~\ref{sec:Main_results}).

The sensitivity of the fronts to noise or a fixed perturbation has also been extensively studied, see e.g.~\cite{Bra83,EckWay94,FifMcL77,Lau85,Sat76,Sat77,Uch78}. In the monostable case, the stability studies lead to a classification of the fronts into two types: ``pulled'' fronts and ``pushed'' fronts~\cite{Rot81,Sto76,Saa03}. A {\it pulled} front is either a critical front $(c^*,U)$ such that the minimal speed $c^*$ satisfies $c^*=2\sqrt{f'(0)},$ or any super-critical front. In the critical case the name pulled front comes from the fact that the front moves at the same speed as the solution of the linearized problem around the unstable state 0, which means that it is being pulled along by its leading edge. This denomination is not so immediate for super-critical fronts. A {\it pushed} front is a critical front $(c^*,U)$  such that  the minimal speed $c^*$ satisfies $c^*>2\sqrt{f'(0)}.$ The speed of propagation of a pushed front is determined not by the behavior of the leading edge of the distribution, but by the whole front. This means that the front is pushed from behind by the nonlinear growth rate in the nonlinear front region itself. A substantial analysis which is not restricted to the monostable case and which relies on the variational structure and the exponential decay of the fronts for more general gradient systems has been carried out in~\cite{LucMurNov04,MurNov08}.

In the present paper, we use a completely different and new approach, which we believe to be simpler and more intuitive, by analyzing the dynamics of the inside structure of the fronts. The results we obtain on the large-time behavior of the components of the fronts in the moving frame and in the whole real line (the precise statements will be given below) shed a new light on and are in keeping with the pulled/pushed terminology in the monostable case as well as with the fact that the bistable or ignition fronts can be viewed as pushed fronts. Even if more general equations or systems could have been considered, we present the results for the one-dimensional equation~\eqref{eq:RD} only, for the sake of simplicity and since this simple one-dimensional situation is sufficient to capture the main interesting properties of the spatial structure of the fronts (however, based on the results of the present paper and on some recent notions of generalized transition waves, we propose in Section~\ref{sec22} some definitions of pulled and pushed transition waves in a more general setting). 

Let us now describe more precisely the model used in this paper. Following the ideas described in~\cite{HalNel08,HalNel09,VlaCav04}, we assume that the fronts are made of several components and we study the behavior of these components. Namely, we consider a traveling wave solution
\begin{equation*}
u(t,x)=U(x-ct)
\end{equation*}
of~\eqref{eq:RD}, where the profile~$U$ satisfies~\eqref{eq:pb_front} and $c$ is the front speed, and we assume that $u$ is initially composed of different groups $(\upsilon_0^i(x))_{i\in I}$ such that, for every $i\in I$,
\begin{equation}\label{eq:u0}
\upsilon_0^i \not\equiv 0\ \hbox{ and }\ 0\leq \upsilon^i_0(x) \leq U(x)  \hbox{ for all } x\in\R,
\end{equation}
where $I$ is a subset of $\N$ and
\[
 u(0,x)=U(x)=\ds\sum_{i\in I}{\upsilon_0^i(x)}\hbox{ for all } x\in\R.
\]
Moreover, all groups $\upsilon^i$ are assumed to share the same characteristics in the sense that they diffuse and grow with the same manner inside the front $u(t,x)$, see~\cite{HalNel08,HalNel09,VlaCav04}. This means that the diffusion coefficient of each group is equal to $1$ and that the per capita growth rate of each group depends only on the entire population and is the same as that of the global front, namely
\begin{equation*}
\ds g(u(t,x)):=\ds \frac{f(u(t,x))}{u(t,x)} \hspace{0.1cm} \hbox{ for all } t\ge0\hbox{ and }x\in\R.
\end{equation*}
In other words, the groups $(\upsilon^i(t,x))_{i\in I}$ satisfy the following equation:
\begin{equation}\label{eq:syst^i}
   \left\{\begin{array}{rcll}
           \ds \partial_t \upsilon^i(t,x) & = & \partial^2_{x}\upsilon^i(t,x) +  g(u(t,x)) \upsilon^i(t,x), & t>0,\ x\in\R,\vspace{3pt} \\
            \upsilon^i(0,x) & = & \upsilon_0^i(x), & x\in\R.
          \end{array}\right.
\end{equation}
Of course it follows from the uniqueness of the solution that
$$\ds u(t,x)=\sum_{i\in I}{\upsilon^i(t,x)} \ \hbox{ for all } t>0 \hbox{ and } x\in\R,$$
which implies that the per capita growth rate $g(u(t,x))= g\big(\sum_{i\in I}{\upsilon^i(t,x)}\big)$ could be viewed as a coupling term in the system~\eqref{eq:syst^i}. The following inequalities also hold from maximum principle
\begin{equation}\label{ineq:u_i}
0<\upsilon^i(t,x)\leq u(t,x) < 1\ \hbox{ for all } t>0,\ x\in\R \hbox{ and  } i\in I.
\end{equation}

Because the components $\upsilon^i$ in \eqref{eq:syst^i} have identical growth and dispersal characteristics, we only have to focus on the behavior of one arbitrarily chosen component $\upsilon^i$ -- we call it $\upsilon$ in the sequel -- to understand the behavior of all the components. This is in sharp contrast with standard competitive systems such as the model of competition between a resident species and an invading species for large open space mentioned in~\cite{ShiKaw97} (see Section~7.2), where usually one of the elements is in some sense stronger than the other one and thus governs the propagation.

Even if the equation~\eqref{eq:RD} is homogeneous and the system~\eqref{eq:syst^i} is linear, one of the main difficulties in our study comes from the fact that a space-time heterogeneity is present in the per capita growth rate $g(u(t,x))$ of each element. It turns out that this heterogeneity does not fulfill any periodicity or monotonicity property. Comparable problems have been studied in~\cite{BerRos07,Ham97_Part1,Ham97_Part2}. In these papers, the authors have considered a reaction term of the form $f(x-ct,\upsilon),$ where the function $\upsilon\mapsto f(x,\upsilon)$ is of the monostable or bistable type for every $x$ and is nonpositive for $\upsilon$ large enough uniformly in $x.$ These properties are not fulfilled here by the function $(t,x,\upsilon)\mapsto g(u(t,x))\,\upsilon$. For instance, if the reaction term $f$ is of type~$(A),$ then $g(u(t,x))$ is always positive. Actually, we prove that the behavior of the groups $\upsilon^i$ mainly depends on the type of $f$, as well as on the initial condition.

The next section is devoted to the statement of our main results. We begin by recalling some asymptotic properties of the solution $U(y)$ of \eqref{eq:pb_front}, as $y\to +\infty.$ Then, the evolution of the density of a group $\upsilon$ solving~\eqref{eq:syst^i} is described in two theorems. Theorem~1 deals with the monostable pulled case and Theorem~2 deals with the monostable pushed case and the bistable and ignition cases. These results show striking differences between the composition of the fronts in the pulled and pushed cases. They lead us to propose new notions of pulled and pushed transition waves in a general setting. The proofs of our results are detailed in Sections~\ref{sec:pulled} and~\ref{sec:pushed}.


\section{Main Results}\label{sec:Main_results}

Let $u(t,x)=U(x-ct)$ be a traveling wave solution of~\eqref{eq:RD} associated to a front $(c,U)$ solving~\eqref{eq:pb_front}, where $f$ is either of type~$(A),$ $(B)$ or~$(C).$ In order to understand the dynamics of a component $\upsilon$ solving~\eqref{eq:u0}-\eqref{eq:syst^i}, inside the traveling wave solution, it is natural to make the following change of variables:
$$\tu(t,x)=\upsilon(t,x+ct) \ \hbox{ for all }t\geq0\hbox{ and } x\in\R.$$
The function $\tu$ corresponds to the solution $\upsilon$ in the moving reference at speed $c$ and it obeys the following equation:
\begin{equation}\label{eq:pb_mv}
\left\{\begin{array}{rcll}
       \partial_t \tu(t,x) & = & \partial^2_{x}\tu(t,x) +  c\partial_x\tu(t,x) + g(U(x))\tu(t,x), & t>0,\ x\in\R,\vspace{3pt} \\
       \tu(0,x) & = & \upsilon_0(x), & x\in\R.
\end{array}\right.
\end{equation}
Thus, the equation \eqref{eq:syst^i} which contains a space-time heterogeneous coefficient reduces to a the reaction-diffusion equation with a spatially heterogeneous coefficient $g(U(x)),$ which only depends on the profile $U$ of the front. As we will see, the leading edge of $U,$ and therefore its asymptotic behavior as $x\to + \infty,$ plays a central role in the dynamics of the solutions of~\eqref{eq:pb_mv}. Before stating our main results, we recall some useful known facts about the asymptotic behavior of the fronts.

\emph{Monostable case (A). }On the one hand, a pulled critical front~$(c^*,U),$ whose speed $c^*$ satisfies $c^*=2\sqrt{f'(0)},$ decays to $0$ at $+\infty$ as follows~\cite{AroWei75,AroWei78}:
\begin{equation}\label{eq:asymp_U_pulled}
\ds U(y)=(Ay+B)\,e^{-\frac{c^*y}{2}}+O\lp e^{-(\frac{c^*}{2}+\delta)y}\rp \ \hbox{ as }y\to+\infty,
\end{equation}
where $\delta>0,$ $A\geq0,$ and $B>0$ if $A=0.$ If $f$ is of the particular KPP type (that is $g(s)=f(s)/s\leq f'(0)$) then $A>0.$ On the other hand, a pushed critical front~$(c^*,U)$, whose speed $c^*$ is such that $c^*>2\sqrt{f'(0)},$ satisfies the following asymptotic property:
\begin{equation}\label{eq:asymp_U_pushed}
\ds U(y)=A\,e^{-\lambda_+(c^*)y}+O\lp e^{-(\lambda_+(c^*)+\delta)y}\rp \ \hbox{ as }y\to+\infty,
\end{equation}
where $\delta>0,$ $A>0$ and $\lambda_+(c^*)=(c^* +\sqrt{(c^*)^2-4f'(0)})/2>c^*/2.$ Thus, the asymptotic behavior of a monostable critical front does depend on its pulled/pushed nature. Lastly, a super-critical front~$(c,U),$ where $c$ satisfies $c>c^*,$ also decays at an exponential rate slower than $c/2$:
\begin{equation}\label{eq:asymp_U_supcritic}
\ds U(y)=A\,e^{-\lambda_-(c)y}+O\lp e^{-(\lambda_-(c)+\delta)y}\rp \ \hbox{ as }y\to+\infty,
\end{equation}
where $\delta>0,$ $A>0$ and $\lambda_-(c)=(c -\sqrt{c^2-4f'(0)})/2<c/2.$

\emph{Bistable case (B). } If follows from~\cite{AroWei75,AroWei78,FifMcL77} that the unique front decays to $0$ at $+\infty$ as follows:
\begin{equation}\label{eq:asymp_U_bi}
       U(y)= A\,e^{-\mu\, y} +O\lp e^{-(\mu+\delta) y}\rp \ \hbox{ as } y\to+\infty,
\end{equation}
where $\delta>0,$ $A>0,$ and $\mu=(c +\sqrt{c^2-4f'(0)})/2>c>c/2.$

\emph{Ignition case (C). } The unique front decays to $0$ at $+\infty$ as follows:
\begin{equation}\label{eq:asymp_U_ing}
       U(y)= A \,e^{-c\, y}  \hbox{ for } y>0 \hbox{ large enough},
\end{equation}
where $A>0$.

We notice that the asymptotic behaviors as $y\to+\infty$ of the fronts in the monostable pushed critical case and the bistable and ignition cases are quite similar. In all cases, the exponential decay rate is faster than $c/2,$ where $c$ is the speed of the front. Let us now state our main results.

\subsection{The inside structure of the fronts}

We first investigate the case where the nonlinearity $f$ is of the monostable type~$(A)$ and $(c,U)$ is a pulled front.

\begin{theo}[Pulled case]\label{ref:theo_mono_pulled}
Let $f$ be of the monostable type~$(A)$, let $(c,U)$ be a pulled front, that is either $c=c^*=2\sqrt{f'(0)}$ or $c>c^*,$ and let $\upsilon$ be the solution of the Cauchy problem~\eqref{eq:syst^i} with the initial condition $\upsilon_0$ satisfying~\eqref{eq:u0} and
\begin{equation}\label{eq:u0_pulled}
    \int_0^{+\infty}{e^{cx}\,\upsilon^2_0(x)\,dx}<+\infty.
\end{equation}
Then
\begin{equation}\label{eq:lim_tilde_ubis}
\limsup_{t\to+\infty}{\Big(\max_{x\geq \alpha\sqrt{t}}{\upsilon(t,x)}\Big)}\to 0\hbox{ as } \alpha\to+\infty.\footnote{Notice that the max in~\eqref{eq:lim_tilde_ubis} is reached from~\eqref{ineq:u_i} and the continuity of~$\upsilon(t,\cdot)$ for all $t>0$.}
\end{equation}
\end{theo}

In other words, any single component $\upsilon$ of the pulled front $u,$ which initially decays faster than the front itself, in the sense of~\eqref{eq:u0_pulled}, cannot follow the propagation of the front. In particular, the formula \eqref{eq:lim_tilde_ubis} implies that
\begin{equation}\label{eq:lim_tilde_u}
  \upsilon(t,x+ct)\to0 \hbox{ uniformly on compacts as } t\to+\infty.
\end{equation}
The conclusion of Theorem~\ref{ref:theo_mono_pulled} holds if  $\upsilon_0$ is of the type $\upsilon_0\equiv U\, \mathds{1}_{(-\infty,a)}$ or more generally if $\upsilon_0$ satisfies~\eqref{eq:u0} and its support is included in $(-\infty,a)$ for some $a\in\R$. This means that the propagation of the traveling wave $u(t,x)=U(x-ct)$ is due to the leading edge of the front.  This characterization agrees with the definition of pulled fronts proposed by Stokes~\cite{Sto76}. It is noteworthy that pulled critical fronts and super-critical fronts share the same inside structure.

Note that \eqref{eq:lim_tilde_ubis} also implies that $\upsilon$ cannot propagate to the right with a positive speed, in the sense that
$$\max_{x\ge\epsilon t}\,\upsilon(t,x)\to0\hbox{ as }t\to+\infty$$
for all $\epsilon>0$. Actually, under some additional assumptions on $\upsilon_0,$ which include the case where $\upsilon_0$ is compactly supported, a stronger uniform convergence result holds:

\begin{prop}\label{ref:prop_pulled}
Under the assumptions of Theorem~$\ref{ref:theo_mono_pulled}$, and if $\upsilon_0$ satisfies the additional condition
\begin{equation}\label{eq:u0_pulled_lim}
  \upsilon_0(x)\to0\hbox{ as }x\to-\infty,\hbox{ or }\upsilon_0\in L^p(\R)\hbox{ for some } p\in[1,+\infty),
\end{equation}
then
\begin{equation}\label{eq:lim_u_pulled}
  \upsilon(t,\cdot)\to 0 \hbox{ uniformly on }\R  \hbox{ as } t\to+\infty.
\end{equation}
\end{prop}

In the pushed case, the dynamics of $\upsilon$ is completely different, as shown by the following result:

\begin{theo}[Pushed case]\label{ref:theo_pushed}
Let $f$ be either of type~$(A)$ with the minimal speed $c^*$ satisfying $c^*>2\sqrt{f'(0)},$ or of type~$(B)$ or $(C).$ Let $(c,U)$ be either the critical front with $c=c^*>2\sqrt{f'(0)}$ in case~$(A)$ or the unique front in case~$(B)$ or~$(C).$ Let $\upsilon$ be the solution of the Cauchy problem~\eqref{eq:syst^i} with the initial condition $\upsilon_0$ satisfying~\eqref{eq:u0}. Then
\begin{equation}\label{eq:sup_left_speed_pushed2}
  \limsup_{t\to+\infty}{\Big( \max_{x\geq\alpha\sqrt{t}}{|\upsilon(t,x)-p(\upsilon_0)U(x-ct)|}\Big)}\to 0\hbox{ as }\alpha\to+\infty,
\end{equation}
where $p(\upsilon_0)\in(0,1]$ is given by
\begin{equation}\label{eq:p}
  p(\upsilon_0)=\ds\frac{\ds\int_\R{\upsilon_0(x)\,U(x)\,e^{cx}\,dx}}{\ds\int_\R{U^2(x)\,e^{cx}\,dx}}.
\end{equation}
Moreover,
\begin{equation}\label{eq:sup_left_speed_pushed}
  \liminf_{t\to+\infty}{\Big(\min_{\alpha\sqrt{t}\leq x\leq x_0+ct}{\upsilon(t,x)}\Big)}>0\hbox{ for all }\alpha\in\R \hbox{ and }x_0\in\R.
\end{equation}
\end{theo}

From~\eqref{eq:asymp_U_pushed},~\eqref{eq:asymp_U_bi} and~\eqref{eq:asymp_U_ing}, $p(\upsilon_0)$ is a well defined positive real number. Theorem~\ref{ref:theo_pushed} is in sharp contrast with Theorem~\ref{ref:theo_mono_pulled}. Indeed, formula \eqref{eq:sup_left_speed_pushed2} in Theorem~\ref{ref:theo_pushed} implies that any small group inside a pushed front is able to follow the traveling wave solution in the sense
\begin{equation}\label{eq:lim_tu_pushed}
  \upsilon(t,x+ct)\to p(\upsilon_0)\,U(x) \hbox{ uniformly on compacts as } t\to+\infty.
\end{equation}
The conclusion~\eqref{eq:lim_tu_pushed} holds even if $\upsilon_0$ is compactly supported. This formula means that an observer who moves with a speed $c$ will see the component $\upsilon$ approach the proportion $p(\upsilon_0)$ of the front $U$. Thus, at large times, the front is made of all its initial components $\upsilon^i_0$ defined in~\eqref{eq:u0}, each one with proportion $p(\upsilon^i_0).$ In other words the front is pushed from the inside. Theorem~2 also shows that the inside structure of the pushed monostable critical fronts and of the bistable and ignition fronts share the same dynamics.

The second formula \eqref{eq:sup_left_speed_pushed} in Theorem~2 shows that the left spreading speed of the group $\upsilon$ inside the front is at least equal to $0$ in the reference frame. More precisely, the group spreads over intervals of the type $(\alpha\sqrt{t},x_0+ct)$ for all $\alpha\in\R$ and $x_0\in\R.$ In fact, the next proposition proves that, if the initial condition $\upsilon_0$ is small at $-\infty,$ then the solution $\upsilon$ spreads to the left with a null speed in the reference frame in the sense that $\upsilon$ is asymptotically small in any interval of the type $(-\infty,\alpha\sqrt{t})$ for $-\alpha>0$ large enough:

\begin{prop}\label{ref:prop_pushed}
Under the assumptions of Theorem~$\ref{ref:theo_pushed}$, if $\upsilon_0$ satisfies the additional assumption~\eqref{eq:u0_pulled_lim}, then
\begin{equation} \label{eq:inf_left_speed}
  \limsup_{t\to+\infty}{\Big( \max_{x\leq\alpha\sqrt{t}}{\upsilon(t,x)}\Big) }\to 0\hbox{ as }\alpha\to-\infty.
\end{equation}
\end{prop}

Notice that without the condition~\eqref{eq:u0_pulled_lim}, the conclusion~\eqref{eq:inf_left_speed} may not hold. Take for instance $\upsilon_0\equiv U,$ then $\upsilon(t,x)=U(x-ct)$ for all $t\ge0$ and $x\in\R,$ and $\sup_{x\le\alpha\sqrt{t}}{\upsilon(t,x)}=1$ for all $\alpha\in\R$ and $t\ge0$.

\begin{remark}{\rm a) One can observe that in the pulled case the function $x\mapsto U(x)e^{cx/2}$ does not belong to $L^2(\R)$, from~\eqref{eq:asymp_U_pulled} and~\eqref{eq:asymp_U_supcritic}. Thus, we can set $p(\upsilon_0)=0$ for any compactly supported initial condition $\upsilon_0$ satisfying~\eqref{eq:u0}. From Theorems~\ref{ref:theo_mono_pulled} and~\ref{ref:theo_pushed}, we can say with this convention that for any monostable reaction term $f$ and any compactly supported $\upsilon_0$ fulfilling~\eqref{eq:u0}, the solution $\upsilon$ of~\eqref{eq:syst^i} is such that $\upsilon(t,x+ct)\to p(\upsilon_0)U(x)$ uniformly on compacts as $t\to+\infty,$ where $p(\upsilon_0)$ is defined by~\eqref{eq:p} if $x\mapsto U(x)e^{cx/2}$ is in $L^2(\R)$ (the pushed case) and $p(\upsilon_0)=0$ otherwise (the pulled case).\par
b) Let us consider the family of reaction terms $(f_a)_{a\ge 0}$ of the monostable type~$(A)$, defined by
$$f_a(u)=u(1-u)(1+au) \hbox{ for all } u\in[0,1] \hbox{ and } a\ge0.$$
The minimal speed $c_a^*$ is given by~\cite{HadRot75}
\[
\ds c_a^*=\left\{\begin{array}{ll}
            2 & \hbox{if } 0\leq a\leq2,\\
            \ds\sqrt{\frac{2}{a}}+\sqrt{\frac{a}{2}} & \hbox{if } a>2.
           \end{array}
\right.
\]
Thus, if $a\in[0,2],$ the critical front $U_a$ associated with $f_a$ is pulled ($c^*_a=2=2\sqrt{f'_a(0)}$) while if $a>2$ the critical front is pushed ($c^*_a>2=2\sqrt{f'_a(0)}$). Up to shift, one can normalize $U_a$ so that $U_a(0)=1/2$ for all $a\ge 0$. A direct computation shows that if $a\geq2$ then the profile of $U_a$ is then given by
\[
\ds U_a(x)=\frac{1}{1+e^{\kappa_a x}}\hbox{ for all }x\in\R,
\]
where  $\kappa_a=\sqrt{a/2}.$ It is easy to check that, if $a>2,$ then the function $x\mapsto U_a(x)e^{c^*_ax/2}$ is in $L^2(\R)$ (a general property shared by all pushed fronts) and $\int_\R{U_a^2(x)\,e^{c^*_a x}\,dx}  \geq  \ds(\kappa_a-\kappa_a^{-1})^{-1}/4 .$ Then, consider a fixed compactly supported initial condition $\upsilon_0$ satisfying~\eqref{eq:u0} and whose support is included in $[-B,B]$ with $B>0.$ Let $p(\upsilon_0,a)$ be defined by~\eqref{eq:p} with $a>2$, $U=U_a$ and $c=c^*_a$. It follows that
\begin{equation*}
0<p(\upsilon_0,a) =  \ds\frac{\ds\int_\R{\upsilon_0(x)\,U_a(x)\,e^{c^*_a x}\,dx}}{\ds\int_\R{U_a^2(x)\,e^{c^*_a x}\,dx}} \leq \ds   4(\kappa_a-\kappa_a^{-1}) \int_{-B}^{B}{U^2_a(x)\,e^{c^*_a x}\,dx} \leq \ds   \frac{8(\kappa_a-\kappa_a^{-1})\cosh{(c_a^*B)}}{c^*_a}.
\end{equation*}
Finally, since $ \kappa_a\to1^+$ and $c^*_a\to2^+$ as $a\to2^+,$ we get that $p(\upsilon_0,a)\to0$ as $a\to2^+.$ Thus, with the convention $p(\upsilon_0)=0$ in the pulled case, this shows the proportion $p(\upsilon_0,a)$ is right-continuous at $a=2,$ which corresponds to the transition between pushed fronts and pulled fronts.}
\end{remark}

\subsection{Notions of pulled and pushed transition waves in a more general setting}\label{sec22}

Our results show that the fronts can be classified in two categories according to the dynamics of their components. This classification agrees with the pulled/pushed terminology introduced by Stokes~\cite{Sto76} in the monostable case and shows that the bistable and ignition fronts have same inside structure as the pushed monostable fronts. This classification also allows us to define the notion of pulled and pushed transition waves in a more general framework. Let us consider the following reaction-dispersion equations:
\begin{equation}\label{eq:pb_general}
           \ds \partial_t u(t,x) = \mathcal{D}(u(t,x)) +  f(t,x,u(t,x)),\ \ t>0,\ x\in\R,
\end{equation}
where $f(t,x,u)$ is assumed to be of class $\C{0,\beta}$ (with $\beta>0$) in $(t, x)$ locally in $u\in[0,+\infty),$ locally Lipschitz-continuous in $u$ uniformly with respect to $(t, x) \in(0,+\infty)\times\R$, $f(\cdot,\cdot,0)=0,$ and $\mathcal{D}$ is a linear operator of dispersion. The classical examples of $\mathcal{D}$ are the homogeneous diffusion operators such as the Laplacian, $\mathcal{D}(u)=D\partial^2_x u$ with $D>0,$ the heterogeneous diffusion operators of the form $\mathcal{D}(u)=\partial_x(a(t,x)\partial_x u)$ where $a(t,x)$ is of class $\C{1,\beta}((0,+\infty)\times\R)$ and uniformly positive, the fractional Laplacian, and the integro-differential operators $\mathcal{D}(u)= J\ast u-u$ where $J\ast u (x)=\int_\R{J(x-y)u(y)dy}$ for all $x\in\R$ and $J$ is a smooth nonnegative kernel of mass $1.$ Before defining the notion of pulled and pushed waves, we recall from~\cite{BerHam11} the definition of transition waves, adapted to the Cauchy problem~\eqref{eq:pb_general}. Let $p^+:(0,+\infty)\times\R\to[0,+\infty)$ be a classical solution of~\eqref{eq:pb_general}. A \emph{transition wave} connecting $p^-=0$ and $p^+$ is a positive solution $u$ of~\eqref{eq:pb_general} such that 1) $u\not\equiv p^\pm$, 2) there exist $n\in\N$ and some disjoint subsets $(\Omega^\pm_t)_{t>0}$ and $(\Gamma_t)_{t>0}=(\{x^1_t,\dots,x^n_t\})_{t>0}$ of $\R$ where $\Gamma_t=\partial\Omega^\pm_t $, $\Omega^-_t\cup\Omega^+_t\cup\Gamma_t=\R$ and $\sup\big\{d(x,\Gamma_t)\,|\ x\in\Omega^-_t\big\}=\sup\big\{d(x,\Gamma_t)\,|\ x\in\Omega^+_t\big\}=+\infty$  for all $t>0$, and 3) for all $\epsilon>0$ there exists $M>0$ such that
$$\ds \hbox{for all } t\in(0,+\infty) \hbox{ and } x\in \overline{\Omega^\pm_t}, \ \Big(d(x,\Gamma_t)\geq M\Big) \Rightarrow \Big(|u(t,x)-p^\pm(t,x)|\leq\epsilon\Big),$$
where $d$ is the classical distance between subsets of $\R$. The wave $u$ for problem~\eqref{eq:pb_general} is also assumed to have a limit $u(0,\cdot)$ at $t=0$ (usually, it is defined for all $t\in\R$, and $f$ and $p^+$ are defined for all $t\in\R$ as well). In the case of Theorems~\ref{ref:theo_mono_pulled} and~\ref{ref:theo_pushed}, the travelling front $u(t,x)=U(x-ct)$ is a transition wave connecting $0$ and $p^+=1,$ the interface $\Gamma_t$ can be reduced to the single point $\Gamma_t=\{x_t\}=\{ct\}$ and the two subsets $\Omega^\pm_t$ can be defined by $\Omega^-_t=(x_t,+\infty)$ and $\Omega^+_t=(-\infty,x_t)$ for all $t>0.$

\begin{defi}[Pulled transition wave]
We say that a transition wave $u$ connecting $0$ and $p^+$ is \emph{pulled} if for any subgroup $\upsilon$ satisfying
\begin{equation}\label{eq:pb_general_comp}
   \left\{\begin{array}{rcll}
           \ds \partial_t \upsilon(t,x) & = & \mathcal{D}(\upsilon(t,x)) +  g(t,x,u(t,x))\upsilon(t,x), & t>0,\ x\in\R, \vspace{3pt}\\
            \upsilon(0,x) & = & \upsilon_0(x), & x\in\R,
          \end{array}\right.
\end{equation}
where $g(t,x,s)=f(t,x,s)/s$ and
\begin{equation}\label{eq:u0_general_comp}
\ds \upsilon_0 \hbox{ is compactly supported}, \ 0\leq \upsilon_0\leq u(0,\cdot)\hbox{ and }\upsilon_0\not\equiv0 ,
\end{equation}
there holds
$$\sup_{d(x,\Gamma_t)\leq M}{\upsilon(t,x)}\to0\hbox{ as } t\to+\infty\hbox{ for all }M\ge 0.$$
\end{defi}

\begin{defi}[Pushed transition wave]
We say that a transition wave $u$ connecting $0$ to $p^+$ is \emph{pushed} if for any subgroup $\upsilon$ satisfying~\eqref{eq:pb_general_comp}-\eqref{eq:u0_general_comp} there exists $M\geq0$ such that
$$\ds  \limsup_{t\to+\infty}{\Big(\sup_{d(x,\Gamma_t)\leq M}{\upsilon(t,x)}\Big)}>0.$$
\end{defi}


\section{The description of pulled fronts}\label{sec:pulled}

We first prove the annihilation of $\upsilon$ in the moving frame, that is property~\eqref{eq:lim_tilde_u}. Then we prove the result~\eqref{eq:lim_tilde_ubis} of Theorem~\ref{ref:theo_mono_pulled} and the result~\eqref{eq:lim_u_pulled} described in Proposition~\ref{ref:prop_pulled} under the additional assumption~\eqref{eq:u0_pulled_lim}. The proof of~\eqref{eq:lim_tilde_u} draws its inspiration from the front stability analyzes in~\cite{Sat76,Sat77,Saa03} and especially from the paper of Eckmann and Wayne~\cite{EckWay94}. It is based on some integral estimates of the ratio $r=\tu/U$ in a suitable weighted space. The proofs of~\eqref{eq:lim_tilde_ubis} and~\eqref{eq:lim_u_pulled} are based on the convergence result~\eqref{eq:lim_tilde_u} and on the maximum principle together with the construction of suitable super-solutions.


\subsection{Local asymptotic extinction in the moving frame: proof of~\eqref{eq:lim_tilde_u}}

Let $f$ be of type~$(A)$ and let $(c,U)$ denote a pulled front satisfying~\eqref{eq:pb_front}, that is $c$ is such that either $c=c^*=2\sqrt{f'(0)}$ or $c>c^*.$ Let  $\upsilon$ be the solution of~\eqref{eq:u0}-\eqref{eq:syst^i} satisfying the condition~\eqref{eq:u0_pulled} and let us set $\tu(t,x)=\upsilon(t,x+ct)$ for all $t\geq0$ and $x\in\R.$ The function $\tu$ solves~\eqref{eq:pb_mv}, while~\eqref{ineq:u_i} implies that
\begin{equation}\label{ineq:tu_U}
0<\tu(t,x)\leq U(x)<1\hbox{ for all } t>0 \hbox{ and }x\in\R.
\end{equation}
Then, let us define the ratio
\begin{equation}\label{eq:r}
\ds r(t,x)=\frac{\tu(t,x)}{U(x)}\hbox{ for all } t\geq0 \hbox{ and }x\in\R.
\end{equation}
The function $r$ is at least of class $\C1$ with respect to $t$ and of class $\C2$ with respect to $x$ in $(0,+\infty)\times\R.$ It satisfies the following Cauchy problem:
\begin{equation}\label{eq:syst^i_r}
\left\{\begin{array}{rcll}
\ds \partial_t{r(t,x)}+\L r(t,x) & = & 0, & t>0,\ x\in\R, \vspace{3pt}\\
\ds r(0,x) & = & \ds\frac{\upsilon_0(x)}{U(x)}, & x\in\R,
\end{array}\right.
\end{equation}
where
\begin{equation*}
\ds\L=-\partial^2_{x}-\psi'(x)\partial_x  \ \hbox{ and }\ \psi(x)=cx+2\ln(U(x))\hbox{ for all } x\in\R.
\end{equation*}

\begin{lem}\label{lem:w_bound}
There exists a constant $k>0$ such that
$$ |\psi'(x)|+|\psi''(x)|\leq k \ \hbox{ for all }x\in\R.   $$
\end{lem}

\noindent\Proof{} The proof uses standard elementary arguments. We just sketch it for the sake of completeness. If we set $q=U$ and $p=U',$ then $ \ds  \psi'=c+2p/q$ and  $\ds \psi''=-2cp/q -2(p/q)^2-2g(q).$ Here we use that $q'=p$ and $p'=-cp-g(q)q.$ Clearly $g(q)$ is bounded. Thus we only need to bound $p/q.$ Proposition~4.4 of~\cite{AroWei78} implies that either $p/q\to-c/2$ at $+\infty$ if $c=c^*=2\sqrt{f'(0)},$ or $p/q\to-\lambda_{-}(c)$ at $+\infty$ if $c>c^*.$ Moreover, since $p/q\to0$ at $-\infty$ and $p/q$ is continuous, we conclude that $p/q$ is bounded, which proves Lemma~\ref{lem:w_bound}.\carre

Let us now define a weight function $\sigma$ as follows:
\begin{equation*}
\sigma(x)= U^2(x)e^{cx}\hbox{ for all } x\in\R.
\end{equation*}
Since $U$ satisfies the asymptotic properties~\eqref{eq:asymp_U_pulled} or~\eqref{eq:asymp_U_supcritic}, one has  $\liminf_{x\to+\infty}{\sigma(x)}>0.$ A direct computation shows that
\begin{equation}\label{eq:sigma}
\sigma'(x)-\psi'(x)\sigma(x)=0\hbox{ for all } x\in\R.
\end{equation}
In order to lighten the proof, we introduce some norms associated to the weight function~$\sigma$:
$$\ds \|w\|_2=\lp\int_{\R}{\sigma(x)\,w(x)^2\,dx}\rp^{1/2}\hbox{ and }\ds \|w\|_\infty=\sup_{x\in\R}{ |\sigma(x)w(x)|},$$
where the supremum is understood as the essential supremum, and we define the standard $L^2(\R)$ and $L^{\infty}(\R)$ norms as follows:
$$\ds |w|_2=\lp\int_{\R}{w(x)^2\,dx}\rp^{1/2}\ \hbox{ and }\ds |w|_\infty=\sup_{x\in\R}{|w(x)|}.$$
Notice that the hypothesis~\eqref{eq:u0_pulled} implies that $r(0,\cdot)$ is in the weighted space
$$L^2_\sigma(\R)=\big\{w\in L^2(\R)\ |\ \|w\|_2<\infty\big\},$$
which is a Hilbert space endowed with the inner product
$$ (w,\widetilde{w})=\int_\R{w(x)\,\widetilde{w}(x)\,\sigma(x)\,dx}\ \hbox{ for all } w,\widetilde{w}\in L^2_\sigma(\R).$$

\begin{lem}\label{lem:coerc_W_Z}
The solution~$r$ of the linear Cauchy problem~\eqref{eq:syst^i_r} satisfies the following properties:
\begin{equation}\label{rrx}
\ds \frac{d}{dt}\lp\frac{1}{2}\|r(t,\cdot)\|^2_2)\rp= -\|\partial_x{r}(t,\cdot)\|^2_2\ \hbox{ for all }t>0
\end{equation}
and, for any constant $K$ satisfying $K\geq|\psi''|_\infty+1,$
\begin{equation*}
\ds \frac{d}{dt}\lp \frac{K}{2}\|r(t,\cdot)\|^2_2+\frac{1}{2}\|\partial_x{r}(t,\cdot)\|^2_2\rp  \leq -\Big( \|\partial_x{r}(t,\cdot)\|^2_2 + \|\partial^2_{x}{r}(t,\cdot)\|^2_2\Big)\hbox{ for all } \  t>0.
\end{equation*}
\end{lem}

\noindent\Proof{} We first set some properties of the operator $\L$ in the Hilbert space $L^2_\sigma(\R).$ We define its domain $D(\L)$ as
$$D(\L)=H^2_\sigma(\R)=\big\{w\in L^2_\sigma(\R)\ |\ w',\,w''\in L^2_\sigma(\R)\big\},$$ where $w'$ and $w''$ are the first- and second-order derivatives of $w$ in the sense of distributions. The domain $D(\L)$ is dense in $L^2_\sigma(\R)$ and since $\C\infty_c(\R)$ is also dense in $H^2_\sigma(\R)$ (with the obvious norm in $H^2_\sigma(\R)$), it follows that
$$ (\L w,w)=\ds -\!\int_\R\!{w''w\sigma}-\!\int_\R\!{\psi'w'w\sigma}=\ds \int_\R\!{w'(w\sigma)'}-\!\int_{\R}\!{\psi'w'w\sigma}=\ds \int_{\R}\!{(w')^2\sigma}+\!\int_\R\!{w'w(\sigma'-\psi'\sigma)}=\|w'\|^2_2\geq0$$
for all  $w\in D(\L)$, that is $\L$ is monotone. Furthermore, $\L$ is maximal in the sense that
\[
\forall\, f\in L^2_\sigma(\R),\  \exists\, w\in D(\L),\ w+\L w=f.
\]
This equation can be solved by approximation (namely, one can first solve the equation $w_n+\L w_n=f$ in $H^2(-n,n)\cap H^1_0(-n,n),$ then show that the sequence $(\|w_n\|_{H^2_\sigma(-n,n)})_{n\in\N}$ is bounded and thus pass to the limit as $n\to+\infty$ to get a solution $w$). Moreover, the operator $\L$ is symmetric since
\[
(\L w,\widetilde{w})=\int_\R{w'(\widetilde{w}\sigma)'}-\int_{\R}{\psi'w'\widetilde{w}\sigma}=\int_\R{w'\widetilde{w}'\sigma}=(w,\L\widetilde{w})
\]
for all $w$, $\widetilde{w}\in D(\L)$. Since $\L$ is maximal, monotone and symmetric, it is thus self-adjoint. Then, since $\upsilon_0$ is in $L^2_\sigma(\R)$, Hille-Yosida Theorem implies that the solution $r$ of the linear Cauchy problem~\eqref{eq:syst^i_r} is such that
\begin{equation}\label{eq:r_dxr_L2}
r\in \C{k}\big((0,+\infty),H^l_\sigma(\R)\big)\cap \C{}  \big([0,+\infty),L^2_\sigma(\R)\big)\hbox{ for all } k,l\in\N,
\end{equation}
where $H^l_\sigma(\R)$ is the set of functions in $L^2_\sigma(\R)$ whose derivatives up to the $l-$th order are in $L^2_\sigma(\R).$

We can now define two additional functions $W$ and $Z.$ First we set,
$$\ds W(t)=\frac{1}{2}\int_{\R}{\sigma(x)r^2(t,x)dx}=\frac{1}{2}\|r(t,\cdot)\|^2_2 \ \hbox{ for all } t\geq0.$$
Let $K$ be a positive constant satisfying $K\geq|\psi''|_\infty+1.$ We define the function $Z$ as follows:
$$\ds Z(t)=KW(t)+\frac{1}{2}\int_{\R}{\sigma(x)\,(\partial_x{r})^2(t,x)\,dx}= \frac{K}{2}\|r(t,\cdot)\|^2_2+\frac{1}{2}\|\partial_x{r}(t,\cdot)\|^2_2\ \hbox{ for all } t>0.$$
The functions $W$ and $Z$ are of the class $\C\infty$ on $(0,+\infty)$ and $W$ is continuous on $[0,+\infty).$ Since $r$ satisfies~\eqref{eq:syst^i_r} and~\eqref{eq:r_dxr_L2} and $\sigma$ obeys~\eqref{eq:sigma}, we get that
\[
W'(t) =  \frac{d}{dt}\Big(\frac{1}{2}\|r(t,\cdot)\|^2_2\Big) =  \frac{1}{2}\frac{d}{dt}\int_{\R}{\sigma(x)\,r^2(t,x)\,dx}  = -\big(\L r(t,\cdot),r(t,\cdot)\big) = -\|\partial_x r(t,\cdot)\|^2_2
\]
for all $t>0$. On the other hand, since $r$ satisfies~\eqref{eq:r_dxr_L2}, there holds
\begin{equation}\label{eq:Z'}
\ds  Z'(t)   =  \ds K W'(t) + \frac{1}{2}\frac{d}{dt}\int_{\R}{\sigma(x)\,(\partial_x r)^2(t,x)\,dx} = KW'(t) + \ds \int_{\R}{\partial_x(\sigma\partial_x r)(t,x)\,\L r(t,x)\,dx}
\end{equation}
for all $t>0$. One can also observe that, for all $t>0,$
\begin{equation*}
\begin{array}{rcl}
\ds-\int_{\R}{\partial_x(\sigma\partial_x r)(t,x)\,\partial^2_x r(t,x)\,dx} & = & - \ds\int_\R{\sigma(x)\,(\partial^2_x r)^2(t,x)\,dx} - \frac{1}{2}\int_\R{\sigma'(x)\,\partial_x((\partial_x r)^2)(t,x)\,dx}\vspace{3pt}\\
& = & - \ds\int_\R{\sigma(x)\,(\partial^2_x r)^2(t,x)\,dx} + \frac{1}{2}\int_\R{\sigma''(x)\,(\partial_x r)^2(t,x)\,dx}
\end{array}
\end{equation*}
and
\begin{equation*}
\begin{array}{l}
\ds-\int_{\R}{\partial_x(\sigma\partial_x r)(t,x)\,\psi'(x)\,\partial_x r(t,x)\,dx}\vspace{3pt}\\
\qquad\qquad\qquad\qquad=\  -\ds\frac{1}{2}\int_\R{\sigma(x)\,\psi'(x)\,\partial_x((\partial_x r)^2)(t,x)\,dx} - \int_\R{\sigma'(x)\,\psi'(x)\,(\partial_x r)^2(t,x)\,dx}\vspace{3pt}\\
\qquad\qquad\qquad\qquad=\ \ds \frac{1}{2}\int_\R{(\sigma \psi')'(x)\,(\partial_x r)^2(t,x)\,dx} -  \int_\R{\sigma'(x)\,\psi'(x)\,(\partial_x r)^2(t,x)\,dx}.
\end{array}
\end{equation*}
Notice that all above integrals exist and all integrations by parts are valid from the density of $\C\infty_c(\R)$ in $H^2_\sigma(\R)$ and from Lemma~\ref{lem:w_bound} and~\eqref{eq:sigma} (in particular, $\sigma'/\sigma$ and $\sigma''/\sigma$ are bounded). Moreover, it follows from~\eqref{eq:sigma} that
\[
\ds \frac{1}{2}(\sigma \psi')'(x)-\sigma'(x)\psi'(x)+\frac{1}{2}\sigma''(x)= \sigma(x)\psi''(x)\leq\sigma(x)|\psi''|_\infty
\]
for all $x\in\R$. Thus, the last integral in~\eqref{eq:Z'} is bounded from above by
\begin{equation*}
\ds\int_{\R}{\partial_x(\sigma\partial_x r)(t,x)\,\L r(t,x)\,dx}\leq  - \ds\int_\R{\sigma(x)\,(\partial^2_x r)^2(t,x)\,dx} + |\psi''|_\infty\int_\R{\sigma(x)\,(\partial_x r)^2(t,x)\,dx}.
\end{equation*}
Then, for all $t>0,$
\begin{equation*}
\begin{array}{rcl}
\ds  Z'(t) & \leq & \ds - \int_\R{\sigma(x)\,(\partial^2_x r)^2(t,x)\,dx} - ( K-|\psi''|_\infty)\int_\R{\sigma(x)\,(\partial_{x}r)^2(t,x)\,dx},\vspace{3pt} \\
         & \leq & \ds - \int_\R{\sigma(x)\,(\partial^2_x r)^2(t,x)\,dx} - \int_\R{\sigma(x)\,(\partial_{x}r)^2(t,x)\,dx},
\end{array}
\end{equation*}
from the choice of $K.$ The proof of Lemma~\ref{lem:coerc_W_Z} is thereby complete.\carre

\noindent{\it{Proof of property~\eqref{eq:lim_tilde_u}.}} We are now ready to state some convergence results for the function $\upsilon$. Note first that $W(t)\geq0$  and $W'(t)\leq0$ for all $t>0.$ Thus $W(t)$ converges to $W_\infty\geq0$ as $t\to+\infty.$ Similarly, $Z(t)$ converges to $Z_\infty\geq0$ as $t\to+\infty$, which implies that $\|\partial_x r(t,\cdot)\|_2$ also converges as $t\to+\infty.$ From~\eqref{rrx} and the convergence of $\|r(t,\cdot)\|_2$ and $\|\partial_xr(t,\cdot)\|_2$ to finite limits, it also follows that
$$\|\partial_xr(t,\cdot)\|_2\to 0\ \hbox{ as }t\to+\infty.$$

On the other hand, for all $t>0$, there holds
\begin{equation}\label{eq:r2_inf}
\ds \|r^2(t,\cdot)\|_\infty\leq\ds\frac{1}{2}\int_{\R}\!{|\partial_x(\sigma r^2)(t,x)|\,dx}\leq\ds \int_{\R}\!{\sigma(x)\,|r(t,x)|\,|\partial_xr(t,x)|\,dx} +\frac{1}{2}\int_{\R}\!{|\sigma'(x)|\,r^2(t,x)\,dx}.
\end{equation}
Using Cauchy-Schwarz inequality, we can see that the first term goes to $0$ as $t\to+\infty$ since $\|\partial_xr(t,\cdot)\|_2$ goes to $0$ and $\|r(t,\cdot)\|_2$ is bounded. If $\sigma'$ were nonnegative on $\R,$ we could drop the modulus in the second term and then, integrating by parts, we would get as above that $\|r^2(t,\cdot)\|_\infty\leq 2\|r(t,\cdot)\|_2\|\partial_x r(t,\cdot)\|_2$ for all $t>0.$

However, the function $\sigma'$ may not be nonnegative on $\R.$ Let us now prove that $\|r^2(t,\cdot)\|_\infty$ still goes to $0$ as $t\to+\infty$ in the general case. First, since $\psi'$ is bounded from Lemma~\ref{lem:w_bound} and $\sigma(x)\sim e^{cx}$ as $x\to-\infty,$ there holds $\sigma'(x)=\psi'(x)\sigma(x)\to0$ as $x\to -\infty.$ Moreover, if $c=c^*=2\sqrt{f'(0)},$ it follows from the asymptotic property~\eqref{eq:asymp_U_pulled} that
\begin{equation*}
\sigma(x)=(Ax+B)^2 +o(1)\hbox{ as }x\to+\infty\ \hbox{ if }\ds U(x)= (A x+B)e^{-c^*x/2} + O\lp e^{-(c^*/2+\delta)x} \rp\hbox{ as }x\to+\infty
\end{equation*}
with $\delta>0,$  $A>0$ and $B\in\R,$ or
\begin{equation*}
\sigma(x)\to B^2>0\hbox{ as }x\to+\infty\ \hbox{ if }\ds U(x)= B e^{-c^*x/2} +O\lp e^{-(c^*/2+\delta)x} \rp\hbox{ as }x\to+\infty,
\end{equation*}
where $\delta>0$ and $B>0.$ On the other hand, if $c>c^*,$ it follows from~\eqref{eq:asymp_U_supcritic},~\eqref{eq:sigma} and the proof of Lemma~\ref{lem:w_bound} that
\begin{equation*}
\sigma'(x)\to+\infty \ \hbox{ as }\ x\to+\infty.
\end{equation*}
Finally, in all cases it is possible to construct a constant $S>0$ and a function $\rho\in\C1(\R)$ such that $\rho'\geq0$ on $\R$ and
\begin{equation*}
S\sigma(x)\leq\rho(x)\leq\sigma(x)\hbox{ for all }x\in\R.
\end{equation*}
The norms associated to $\rho$ are equivalent to those defined by $\sigma.$ Then, denoting
\[
\ds \|w\|_{\rho,2}=\lp\int_{\R}{\rho(x)w^2(x)dx}\rp^{1/2} \ \hbox{ and }\ \|w\|_{\rho,\infty}=\sup_{x\in\R}{ |\rho(x)w(x)|},
\]
and applying equation~\eqref{eq:r2_inf} to these norms, one infers that
\[
S\|r^2(t,\cdot)\|_\infty\leq \|r^2(t,\cdot)\|_{\rho,\infty}\leq 2\|r(t,\cdot)\|_{\rho,2}\|\partial_xr(t,\cdot)\|_{\rho,2}\leq2\|r(t,\cdot)\|_{2}\|\partial_xr(t,\cdot)\|_{2}
\]
for all $t>0$. Since $\|r(t,\cdot)\|_2$ is bounded and $\|\partial_xr(t,\cdot)\|_2\to0$ as $t\to+\infty,$ it follows that
\begin{equation}\label{eq:limsup_r2}
\lim_{t\to+\infty}{\Big(\sup_{x\in\R}{\big(\sigma(x)r^2(t,x)\big)}\Big)}=0.
\end{equation}
Moreover,~\eqref{eq:r} implies that
\begin{equation*}
0\leq \tu(t,x)=\big(U^2(x)r^2(t,x)\big)^{1/2}=\big(\sigma(x)r^2(t,x)\big)^{1/2}e^{-cx/2}
\end{equation*}
for all $t>0$ and $x\in\R$. Then, for any compact set $\mathcal{K},$ one has
\begin{equation*}
\max_{x\in\mathcal{K}}\,{\upsilon(t,x+ct)}=\ds \max_{x\in\mathcal{K}}\,{\tu(t,x)}\leq \Big(\max_{x\in\mathcal{K}}\,e^{-cx/2}\Big)\times\Big(\sup_{x\in\R}{\big(\sigma(x)r^2(t,x)\big)^{1/2}}\Big).
\end{equation*}
Finally, equation~\eqref{eq:limsup_r2} implies that $\tu$ converge to $0$ uniformly on compacts as $t\to+\infty,$ which yields~\eqref{eq:lim_tilde_u}. One can also say that
\begin{equation}\label{eq:lim_u_pulled_mv}
\lim_{t\to+\infty}{\Big(\max_{x\geq A}\,{\upsilon(t,x+ct)}\Big)}=\lim_{t\to+\infty}{\Big(\max_{x\geq A+ct}{\upsilon(t,x)}\Big)}=0.
\end{equation}
for all $A\in\R$, where we recall that the maxima in~\eqref{eq:lim_u_pulled_mv} are reached from~\eqref{ineq:u_i} and the continuity of $\upsilon(t,\cdot)$ for all $t>0$.\hfill$\Box$


\subsection{Extinction in $(\alpha\sqrt{t},+\infty)$ and in $\R$: proofs of Theorem~\ref{ref:theo_mono_pulled} and Proposition~\ref{ref:prop_pulled}}

Before completing the proof of Theorem~\ref{ref:theo_mono_pulled} and Proposition~\ref{ref:prop_pulled}, let us first state two auxiliary lemmas. They provide some uniform estimates of $\upsilon$ in intervals of the type $(\alpha\sqrt{t},A+ct)$ or the whole real line $\R$ when bounds for $\upsilon(t,\cdot)$ are known at the positions $A+ct$, in the intervals $(A+ct,+\infty)$ and/or at $-\infty$. These two lemmas will be used in all cases~$(A)$, $(B)$ and~$(C).$

\begin{lem}\label{lem:bounds_u}
Let $f$ be of type~$(A),$ $(B)$ or~$(C),$ let $(c,U)$ be a front satisfying~\eqref{eq:pb_front} and let $\upsilon$ solve~\eqref{eq:syst^i} with $\upsilon_0$ satisfying~\eqref{eq:u0}. Let $\mu\in[0,1]$ and $A_0\in\R$ be such that $\limsup_{t\to+\infty}\upsilon(t,A+ct)\le\mu$ $($resp. $\liminf_{t\to+\infty}\upsilon(t,A+ct)\ge\mu)$ for all $A<A_0$. Then, for all $\epsilon>0,$ there exist $\alpha_0>0$ and $A<A_0$ such that
\begin{equation}\label{eq:bound_u_sup}
\limsup_{t\to+\infty}{\Big(\max_{\alpha\sqrt{t}\leq x \leq A+ct}{\upsilon(t,x)}\Big)}\leq\mu+\epsilon \ \hbox{ for all  }\alpha\geq\alpha_0,
\end{equation}
resp.
$$\liminf_{t\to+\infty}{\Big(\min_{\alpha\sqrt{t}\leq x \leq A+ct}{\upsilon(t,x)}\Big)}\geq\mu-\epsilon \ \hbox{ for all  }\alpha\geq\alpha_0.$$
\end{lem}

\begin{lem}\label{lem:sup_tu<lambda}
Let $f$ be of type~$(A),$ $(B)$ or~$(C),$ let $(c,U)$ be a front satisfying~\eqref{eq:pb_front} and let $\upsilon$ solve~\eqref{eq:syst^i} with $\upsilon_0$ satisfying~\eqref{eq:u0}. Let $\lambda$ and $\mu\in[0,1]$ be such that either $\limsup_{x\to-\infty}{\upsilon_0(x)}\le\lambda$ or $\max{(\upsilon_0-\lambda,0)}\in L^p(\R)$ for some $p\in[1,+\infty)$, and $\limsup_{t\to+\infty}\big(\max_{x\ge A+ct}\upsilon(t,x)\big)\le\mu$ for all $A\in\R$. Then
\begin{equation}\label{eq:sup_tu<lambda}
\ds \limsup_{t\to+\infty}{\Big(\sup_{x\in\R}{\upsilon(t,x)}\Big)}\leq\max{(\lambda,\mu)}.
\end{equation}
\end{lem}

The proofs of Lemmas~\ref{lem:bounds_u} and~\ref{lem:sup_tu<lambda} are postponed at the end of this section.\hfill\break

\noindent{\it{End of the proof of Theorem~$\ref{ref:theo_mono_pulled}$.}} Let $\upsilon$ be the solution of~\eqref{eq:syst^i} with $\upsilon_0$ satisfying~\eqref{eq:u0} and~\eqref{eq:u0_pulled}. To get~\eqref{eq:lim_tilde_ubis}, pick any $\epsilon>0$ and observe that property~\eqref{eq:lim_tilde_u} and Lemma~\ref{lem:bounds_u} with $\mu=0$ (and an arbitrary $A_0$) yield the existence of $\alpha_0>0$ and $A\in\R$ such that
\begin{equation*}
\limsup_{t\to+\infty}{\Big(\max_{\alpha\sqrt{t}\leq x\leq A+ct}{\upsilon(t,x)}\Big)}\leq\epsilon \hbox{ for all } \alpha\geq\alpha_0.
\end{equation*}
Property~\eqref{eq:lim_tilde_ubis} follows then from~\eqref{eq:lim_u_pulled_mv}. This proves Theorem~$\ref{ref:theo_mono_pulled}.$\carre

\noindent{\it{End of the proof of  Proposition~$\ref{ref:prop_pulled}$.}} We make the additional assumption~\eqref{eq:u0_pulled_lim}. Notice that the assumptions of Lemma~\ref{lem:sup_tu<lambda} are fulfilled with $\lambda=\mu=0$, from~\eqref{eq:u0_pulled_lim} and~\eqref{eq:lim_u_pulled_mv}. It follows that the inequality~\eqref{eq:sup_tu<lambda} holds with $\lambda=\mu=0,$ which implies that $\upsilon(t,x)\to0$ uniformly on $\R$ as $t\to+\infty$. The proof of Proposition~\ref{ref:prop_pulled} is thereby complete.\carre

The proofs of Lemmas~\ref{lem:bounds_u} and~\ref{lem:sup_tu<lambda} are based on the construction of explicit sub- or super-solutions of~\eqref{eq:syst^i} in suitable domains in the $(t,x)$ coordinates, and on the fact that the coefficient $g(u(t,x))=g(U(x-ct))$ in~\eqref{eq:syst^i} is exponentially small when $x-ct\to-\infty$.\hfill\break

\noindent{\it{Proof of Lemma~$\ref{lem:bounds_u}$.}} Let us first consider the case of the upper bounds. Let $\mu\in[0,1]$ and $A_0\in\R$ be as in the statement and pick any $\epsilon\in(0,1).$ Let $j_\epsilon$ be the positive function defined on $(-\infty,0)$ by
\begin{equation}\label{eq:j_epsilon}
j_\epsilon(y)=\epsilon\lp 1-\frac{1}{1-y} \rp\ \hbox{ for all } y<0.
\end{equation}
The functions $j_\epsilon'$ and $j_\epsilon''$ are negative on $(-\infty,0)$ and $-j_\epsilon''(y)-cj_\epsilon'(y)\sim\epsilon c/y^2$ as $y\to-\infty.$ We recall the asymptotic behavior of $U$ at $-\infty$:
\begin{equation}\label{eq:asymp_U_-infty}
       U(y)= 1 - Be^{\nu y}+O\lp e^{(\nu+\delta) y}\rp \  \hbox{ as } y\to-\infty,
\end{equation}
where $\delta>0,$ $B>0$ and $\nu=(-c +\sqrt{c^2-4f'(1)})/2>0.$ This property~\eqref{eq:asymp_U_-infty} and the negativity of $f'(1)$ yield the existence of a real number $A<\min(A_0,0)$ such that
\begin{equation}\label{ineq:j_A}
\ds 0\leq g(U(y))\leq \frac{-j_\epsilon''(y)-cj_\epsilon'(y)}{2+\epsilon}\ \hbox{ for all } y\in(-\infty,A].
\end{equation}
From the assumptions made in Lemma~\ref{lem:bounds_u}, there is $t_0>0$ such that
$$\upsilon(t,A+ct)\leq\mu+\epsilon\hbox{ for all }t\geq t_0.$$

Now, let us define the function $\overline{\upsilon}$ by
$$\overline{\upsilon}(t,x)=h(t-t_0,x-A-ct_0)+j_\epsilon(x-ct)\hbox{ for all  } t\geq t_0\hbox{ and } x\in[A+ct_0,A+ct],$$
where $h$ solves the heat equation
\begin{equation*}
\left\{\begin{array}{rcll}
    \ds \partial_t h(t,x) & = & \partial_x^2 h(t,x), & t>0,\ x\in\R, \vspace{3pt}\\
    \ds h(0,x) & = & 2\,\mathds{1}_{(-\infty,0)}(x)+(\mu+\epsilon)\mathds{1}_{(0,+\infty)}(x), & x\in\R.
\end{array}\right.
\end{equation*}
Let us check that $\overline{\upsilon}$ is a supersolution of~\eqref{eq:syst^i} in the domain $t\geq t_0$ and $x\in[A+ct_0,A+ct].$ Firstly, observe that
\begin{equation}\label{bounds}
\mu+\epsilon\le\overline{\upsilon}(t,x)\le 2+\epsilon
\end{equation}
for all $t\ge t_0$ and $x\in[A+ct_0,A+ct]$ since $\mu+\epsilon<2$ and $0<j_{\epsilon}<\epsilon$ on $(-\infty,0)$. It follows from~\eqref{ineq:tu_U} that
$$\upsilon(t,A+ct_0)<1\le\frac{2+\mu+\epsilon}{2}=h(t-t_0,0)\le\overline{\upsilon}(t,A+ct_0)\hbox{ for all }t>t_0.$$
On the other hand,
\begin{equation*}
\upsilon(t,A+ct)\le\mu+\epsilon\le h(t-t_0,c(t-t_0))\le\overline{\upsilon}(t,A+ct)\hbox{ for all } t\geq t_0.
\end{equation*}
Lastly, from~\eqref{ineq:j_A} and~\eqref{bounds}, the function $\overline{\upsilon}$ satisfies, for all $t>t_0$ and $x\in(A+ct_0,A+ct)$,
\begin{equation}\label{eq:u_ou}
\begin{array}{l}
\partial_t \overline{\upsilon}(t,x) - \partial_x^2 \overline{\upsilon}(t,x)  -g(U(x-ct))\overline{\upsilon}(t,x) \vspace{3pt}\\
\qquad\qquad\qquad=\ -j_\epsilon''(x-ct)-cj_\epsilon'(x-ct)-g(U(x-ct))\overline{\upsilon}(t,x) \vspace{3pt}\\
\qquad\qquad\qquad\geq\ -g\big(U(x-ct)\big)(2+\epsilon) -j_\epsilon''(x-ct)-cj_\epsilon'(x-ct) \vspace{3pt}\\
\qquad\qquad\qquad\geq\ 0.
\end{array}
\end{equation}
The maximum principle applied to~\eqref{eq:syst^i} implies that
$$\upsilon(t,x)\leq \overline{\upsilon}(t,x)\hbox{ for all }t\geq t_0\hbox{ and }x\in[A+ct_0,A+ct].$$

Fix $\alpha_0>0$ so that
\begin{equation*}
\ds \frac{2-\mu-\epsilon}{\sqrt{\pi}}\ds \int^{-\alpha_0/2}_{-\infty}{\ds e^{-z^2}dz}\leq \epsilon.
\end{equation*}
Let $\alpha\geq\alpha_0$ and $t_1>t_0$ be such that $A+ct_0<\alpha\sqrt{t}<A+ct$ for all $t\geq t_1.$ Since $h(t,\cdot)$ is decreasing for all $t>0$, there holds, for all $t\geq t_1$,
\begin{equation*}
\begin{array}{rcl}
\ds \max_{\alpha\sqrt{t}\leq x \leq A+ct}{\upsilon(t,x)} & \leq & \ds \max_{\alpha\sqrt{t}\leq x \leq A+ct}{\overline{\upsilon}(t,x)} \\[0.3cm]
                                                 & = & \ds \max_{\alpha\sqrt{t}\leq x \leq A+ct}{\big(h(t-t_0,x-A-ct_0)+j_\epsilon(x-ct)\big)} \\[0.3cm]
                                                 & \leq & \ds h(t-t_0,\alpha\sqrt{t}-A-ct_0)+\epsilon \\
                                                 & \leq & \ds \mu+2\epsilon +\frac{2-\mu-\epsilon}{\sqrt{\pi}}\ds \int^{(-\alpha\sqrt{t}+A+ct_0)/\sqrt{4(t-t_0)}}_{-\infty}{\ds e^{-y^2}dy}.
\end{array}
\end{equation*}
Hence, for all $\alpha\ge\alpha_0$,
\begin{equation*}
\limsup_{t\to+\infty}{\Big(\max_{\alpha\sqrt{t}\leq x \leq A+ct}{\upsilon(t,x)}\Big)}\leq \mu+2\epsilon +\frac{2-\mu-\epsilon}{\sqrt{\pi}}\ds \int^{-\alpha/2}_{-\infty}{\ds e^{-y^2}dy} \leq \mu+3\epsilon.
\end{equation*}
Since $\epsilon>0$ is arbitrarily small, the conclusion~\eqref{eq:bound_u_sup} follows.

As far as the lower bounds are concerned, with the same type of arguments as above and since $g(U(y))$ is nonnegative near $-\infty,$ one can show that for any $\epsilon\in(0,1),$ there exist $A<A_0$ and $t_0>0$ such that
$$\upsilon(t,x)\geq h(t-t_0,x-A-ct_0)\hbox{ for all }t\geq t_0\hbox{ and }x\in[A+ct_0,A+ct],$$
where $h$ solves the heat equation in $\R$ with initial condition $h(0,\cdot)=-\mathds{1}_{(-\infty,0)}+(\mu-\epsilon)\mathds{1}_{(0,+\infty)}$. It follows that
\begin{equation*}
\liminf_{t\to+\infty}{\Big(\min_{\alpha\sqrt{t}\leq x \leq A+ct}{\upsilon(t,x)}\Big)}\geq \mu-2\epsilon \hbox{ for } \alpha>0 \hbox{ large enough.}
\end{equation*}
The proof of Lemma~\ref{lem:bounds_u} is thereby complete.\carre

\noindent{\it{Proof of Lemma~$\ref{lem:sup_tu<lambda}$.}} Let $\lambda\in[0,1]$ and $\mu\in[0,1]$ be given as in the statement and pick any $\epsilon>0.$ Let $j_\epsilon$ be the function defined as in~\eqref{eq:j_epsilon} and let $A<0$ be such that
\begin{equation}\label{ineq:j}
\ds 0\leq g(U(y))\leq \frac{-j_\epsilon''(y)-cj_\epsilon'(y)}{\max(\max(\lambda,\mu)+\epsilon,1)+\epsilon}\hbox{ for all } y\in(-\infty,A].
\end{equation}
The assumptions made in Lemma~\ref{lem:sup_tu<lambda} yield the existence of $t_0>0$ such that
\begin{equation}\label{eq:tu_A_epsilon}
\upsilon(t,x)\leq\mu +\epsilon\hbox{ for all } t\geq t_0\hbox{ and } x\in[A+ct,+\infty).
\end{equation}
If $\upsilon_0$ satisfies $\limsup_{x\to-\infty}{\upsilon_0(x)}\leq\lambda,$ then the comparison with the heat equation and the equality $g(1)=0$ imply that $\limsup_{x\to-\infty}{\upsilon(t,x)}\leq\lambda$ for all $t>0.$ On the other hand, if $w$ solves~\eqref{eq:syst^i} with an initial condition $w_0$ in $L^p(\R)$ for some $p\in[1,+\infty),$ then heat kernel estimates and the boundedness of $g(U(x-ct))$  imply that $w(t,\cdot)$ is also in $L^p(\R)$ for all $t>0$, while $w(t,\cdot)$ is uniformly continuous on $\R$ from standard parabolic estimates. Finally,  $w(t,x)\to0$ as $x\to-\infty$ for all $t>0.$ Now, if $\max{(\upsilon_0-\lambda,0)}\in L^p(\R)$ for some $p\in[1,+\infty),$ then by writing $\upsilon_0\leq\lambda+\max{(\upsilon_0-\lambda,0)},$ the previous arguments and the linearity of~\eqref{eq:syst^i} imply that $\limsup_{x\to-\infty}{\upsilon(t,x)}\leq\lambda$ for all $t>0.$ In any case, at time $t=t_0,$ there exists $B\leq A$ such that
\begin{equation}\label{eq:tu_B_epsilon}
\upsilon(t_0,x)\leq \lambda +\epsilon\hbox{ for all }  x\in(-\infty,B+ct_0].
\end{equation}

Now, let us define the function $\overline{\upsilon}$ by
\begin{equation*}
\ds \overline{\upsilon}(t,x)=h(t-t_0,x)+j_\epsilon(x-ct)\hbox{ for all } t\geq t_0 \hbox{ and } x\in(-\infty,A+ct],
\end{equation*}
where $h$ solves the heat equation
\begin{equation}\label{heat}
\left\{\begin{array}{rcll}
    \ds \partial_t h(t,x) & = & \partial_x^2 h(t,x), & t>0,\ x\in\R,\vspace{3pt} \\
    \ds h(0,x) & = & \max(\max(\lambda,\mu)+\epsilon,\upsilon(t_0,x)), & x\in\R.
\end{array}\right.
\end{equation}
Let us check that $\overline{\upsilon}$ is a supersolution of~\eqref{eq:syst^i} in the domain $t\geq t_0$ and $x\leq A+ct.$ Firstly, observe that
\begin{equation}\label{boundsbis}
\max{(\lambda,\mu)}+\epsilon\le\overline{\upsilon}(t,x)\le\max(\max{(\lambda,\mu)}+\epsilon,1)+\epsilon
\end{equation}
for all $t\geq t_0$ and $x\leq A+ct$, from~\eqref{ineq:tu_U},~\eqref{eq:j_epsilon} and the maximum principle applied to the heat equation~\eqref{heat}. Then, from equation~\eqref{eq:tu_A_epsilon}, we get that
$$\upsilon(t,A+ct)\leq\mu+\epsilon\leq\overline{\upsilon}(t,A+ct) \hbox{ for all } t\geq t_0.$$
Moreover, by definition of $\overline{\upsilon}(t_0,\cdot),$ there holds
\[
\upsilon(t_0,x)\leq \max(\max{(\lambda,\mu)}+\epsilon,\upsilon(t_0,x)) \leq \overline{\upsilon}(t_0,x)\hbox{ for all }x\in(-\infty,A+ct_0].
\]
Finally, from~\eqref{ineq:j} and~\eqref{boundsbis}, $\overline{\upsilon}$ satisfies the following inequality, for all $t>t_0$ and $x\in(-\infty,A+ct)$,
\begin{equation*}
\begin{array}{l}
\partial_t \overline{\upsilon}(t,x) - \partial_x^2 \overline{\upsilon}(t,x)  -g(U(x-ct))\overline{\upsilon}(t,x)\vspace{3pt}\\
\qquad\qquad\qquad=\ -j_\epsilon''(x-ct)-cj_\epsilon'(x-ct)-g(U(x-ct))\overline{\upsilon}(t,x)\vspace{3pt} \\
\qquad\qquad\qquad\geq\ -g(U(x-ct))\big(\max(\max(\lambda,\mu)+\epsilon,1)+\epsilon\big) -j_\epsilon''(x-ct)-cj_\epsilon'(x-ct)\vspace{3pt} \\
\qquad\qquad\qquad\geq\ 0.
\end{array}
\end{equation*}
The maximum principle applied to~\eqref{eq:syst^i} implies that
\begin{equation}\label{eq:tu<ou}
\upsilon(t,x)\leq \overline{\upsilon}(t,x)\hbox{ for all } t\geq t_0\hbox{ and } x\in(-\infty,A+ct].
\end{equation}

Next, we claim that
$$\overline{\upsilon}(t,x)-j_\epsilon(x-ct)\to\max{(\lambda,\mu)}+\epsilon\hbox{ uniformly on }(-\infty,A+ct]\hbox{ as }t\to+\infty.$$
Indeed, since $\upsilon(t_0,\cdot)$ satisfies~\eqref{eq:tu_A_epsilon} and~\eqref{eq:tu_B_epsilon}, the initial condition $h(0,\cdot)$ of $h$ is the sum of the constant $\max{(\lambda,\mu)}+\epsilon$ and a nonnegative compactly supported continuous function. By linearity and standard properties of the heat equation on $\R,$ it follows that $h(t,\cdot)\to\max{(\lambda,\mu)}+\epsilon$ uniformly on $\R$ as $t\to+\infty.$ Hence, $\overline{\upsilon}(t,x)-j_\epsilon(x-ct)\to\max{(\lambda,\mu)}+\epsilon$ uniformly on $(-\infty,A+ct]$ as $t\to+\infty$. From~\eqref{eq:tu<ou}, we get that
\[
\limsup_{t\to+\infty}{\Big(\sup_{x\in(-\infty,A+ct]}{\upsilon(t,x)}\Big)}\leq \max{(\lambda,\mu)}+\epsilon +\sup_{y\in(-\infty,A]}j_\epsilon(y)=\max{(\lambda,\mu)}+2\epsilon.
\]
Then, from equation~\eqref{eq:tu_A_epsilon} we get that
\[
\limsup_{t\to+\infty}{\Big(\sup_{x\in\R}{\upsilon(t,x)}\Big)} \leq \max{(\lambda,\mu)} +2\epsilon.
\]
Since $\epsilon>0$ is arbitrary,  the conclusion~\eqref{eq:sup_tu<lambda} follows and the proof of Lemma~\ref{lem:sup_tu<lambda} is complete.\hfill$\Box$

\begin{remark}\label{rem:1}{\rm 
If the initial condition $\upsilon_0$ of~\eqref{eq:syst^i} is such that $\upsilon_0(x)\to\lambda\in[0,1]$ as $x\to-\infty,$ then, as already noticed, $\upsilon(t,x)\to\lambda$ as $x\to-\infty$ for all $t>0,$ since $g(1)=0.$ In particular, $\sup_{\R}{\upsilon(t,\cdot)}\geq\lambda$ for all $t>0.$ Therefore, if $\upsilon_0$ satisfies~\eqref{eq:u0},~\eqref{eq:u0_pulled} and $\lim_{x\to-\infty}{\upsilon_0(x)}=\lambda\in[0,1],$ then the proof of Lemma~$\ref{lem:sup_tu<lambda}$ shows that $\sup_\R{\upsilon(t,\cdot)}\to\lambda$ as $t\to+\infty.$}
\end{remark}


\section{The description of pushed fronts}\label{sec:pushed}

This section is devoted to the proofs of Theorem~\ref{ref:theo_pushed} and Proposition~\ref{ref:prop_pushed}. We begin by proving formula~\eqref{eq:lim_tu_pushed}. The proof of this formula draws its inspiration from the front stability analysis~\cite{Sat76,Sat77} and especially from the lecture notes of Gallay~\cite{Gal05}. It is based on some properties of the self-adjoint Schr\"{o}dinger operator $\L$ defined by
\begin{equation}\label{eq:def_L}
\ds \L=-\partial^2_{x}+\lp \frac{{c}^2}{4}-g(U(x))\rp
\end{equation}
with domain $H^2(\R),$ where $g(s)=f(s)/s,$ $f$ satisfies the hypotheses of Theorem~\ref{ref:theo_pushed}, and $(c,U)$ is either the pushed critical front in case~$(A)$ when $c=c^*>2\sqrt{f'(0)},$ or the unique front satisfying~\eqref{eq:pb_front} in case~$(B)$ or $(C).$ The properties of the semigroup generated by $-\L$ play an essential role in the large time behavior of the solution $\tu$ of the Cauchy problem~\eqref{eq:pb_mv}. Indeed, the function $\upsilon^*$ defined by
\begin{equation}\label{eq:def_su*}
\su(t,x)=e^{cx/2}\,\tu(t,x)\hbox{ for all } t\geq0 \hbox{ and }x\in\R
\end{equation}
satisfies the Cauchy problem
\begin{equation}\label{eq:pb_su*}
\left\{\begin{array}{rcll}
       \partial_t \su(t,x) + \L \su(t,x) & = & 0, & t>0,\ x\in\R, \vspace{3pt}\\
       \su(0,x) & = & \upsilon_0(x)\,e^{cx/2}, & x\in\R.
\end{array}\right.
\end{equation}
The main spectral properties of $\L$ are stated in Section~\ref{sec41}. Then, Section~\ref{sec42} is devoted to the proof of formula~\eqref{eq:lim_tu_pushed}. The proofs of Theorem~\ref{ref:theo_pushed} and Proposition~\ref{ref:prop_pushed} are given in Section~\ref{sec:spread}.


\subsection{Preliminary lemmas}\label{sec41}

Let $X_{c/2}$ be the weighted space defined by
\begin{equation}\label{eq:def_X}
\ds X_{c/2}=\left\{ w\in L^2_{loc}(\R)\,|\,\int_{\R}{w^2(x)\,e^{cx}\,dx}<\infty\right\}.
\end{equation}

\begin{lem}\label{lem:L1}
If the front $U$ solving~\eqref{eq:pb_front} with $c>0$ belongs to $X_{c/2},$ then the operator $\L$ defined by~\eqref{eq:def_L} satisfies the following properties:
\begin{itemize}
\item[i)] The essential spectrum $\sigma_e(\L)$ of $\L$ is equal to $\big[c^2/4-\max{(f'(0),0)},+\infty\big).$
\item[ii)] The point spectrum of $\L$ is included in $\big[0,c^2/4-\max{(f'(0),0)}\big).$ Moreover, $\lambda=0$ is the smallest eigenvalue of $\L$ and the function $x\mapsto\phi(x)=U(x)\,e^{cx/2}$ spans the kernel of~$\L$.
\item[iii)] The following spectral decomposition of $L^2(\R)$ holds:
\begin{equation}\label{eq:decomp_L2}
L^2(\R)=\im(P)\oplus\ker(P),
\end{equation}
where the operator $P:L^2(\R)\to L^2(\R)$ is the spectral projection onto the kernel of $\L,$ that is $P(w)=\big(\int_{\R}\varphi w\big)\,\varphi$ for all $w\in L^2(\R)$ with $\varphi=\phi/\|\phi\|_{L^2(\R)}$.
\end{itemize}
\end{lem}

\noindent\Proof{} It uses standard results and it is just sketched here for the sake of completeness.

{\it{i)}} The coefficients of the operator $\L$ are not constant but converge exponentially to two limits as $x\to\pm\infty.$ It follows that $\L$ is a relatively compact perturbation of the operator $\L_0$ defined by $\L_0 = -\partial^2_{x} +(c^2/4-g_\infty(x))$, where
\begin{equation*}
g_\infty(x)=\left\{\begin{array}{ll}
                   g(0)=f'(0) & \hbox{if }x<0,\\
                   g(1)=0     & \hbox{if }x>0.
                  \end{array}\right.
\end{equation*}
Then, Theorem~A.2 in~\cite{Hen81} implies that the essential spectrum $\sigma_e(\L)$ of the operator $\L$ is equal to the spectrum $\sigma(\L_0)$ of the operator $\L_0.$ Since $\L_0$ is self-adjoint in $L^2(\R),$ we get
\begin{equation*}
\sigma_e(\L)=\sigma(\L_0)=\big[c^2/4-\max{(f'(0),0)},+\infty\big).
\end{equation*}

{\it{ii)}} The operator $\L$ is a self-adjoint operator in $L^2(\R),$ so the eigenvalues of $\L$ are in $\R.$ Moreover, since $(c,U)$ satisfies equation~\eqref{eq:pb_front} and $U$ belongs to $X_{c/2},$ one has necessarily that $c^2>4f'(0)$ (whatever the sign of $f'(0)$ be) and the function $x\mapsto\phi(x)=U(x)e^{cx/2}$ is in $L^2(\R)$ (and then in $H^2(\R)$ by adapting the arguments used the proof of Lemma~\ref{lem:w_bound}). Furthermore, $\phi$ is an eigenvector of $\L$ associated to the eigenvalue $\lambda=0,$ that is $\L\phi=0$, and the eigenvalue is simple, from elementary arguments based on the exponential behavior at $\pm\infty$. On the other hand, since $\phi$ is positive, Sturm-Liouville theory implies that $\lambda=0$ is the lowest value of the spectrum of $\L$. Together with {\it{i}}), we finally get that the point spectrum of $\L$ is a discrete subset of the interval $\big[0,c^2/4-\max(f'(0),0)\big).$

{\it{iii)}} The function $\varphi=\phi/\|\phi\|_{L^2(\R)}$ is a normalized eigenvector of $\L$ associated to the eigenvalue $0.$ Since $\L$ is self-adjoint, the operator $P:L^2(\R)\to L^2(\R)$ defined as in Lemma~\ref{lem:L1} is the spectral projection onto the kernel of $\L.$ Then, the spectral decomposition~\eqref{eq:decomp_L2} holds, where $ \im(P)=\big\{\beta\varphi,\ \beta\in\R\big\}$ and $\ker(P)=\big\{ w\in L^2(\R)\,|\, \int_\R{\varphi w}=0 \big\}.$\carre

Let us come back to the Cauchy problem $\partial_t w+\L w=0.$ From Lemma~\ref{lem:L1}, the semigroup~$\big(e^{-t\L}\big)_{t\geq0}$ generated by~$-\L$ satisfies the following properties:

\begin{lem}\label{lem:L2}
If the front $U$ solving~\eqref{eq:pb_front} with $c>0$ belongs to $X_{c/2}$, then there exist two constants $C>0$ and $\eta>0$ such that
\begin{equation}\label{eq:lim_kerP}
|e^{-t\L}w|_\infty\leq C e^ {-\eta t}|w|_\infty \hbox{ for all }t\geq0 \hbox{ and } w\in\ker(P)\cap L^\infty(\R).
\end{equation}
\end{lem}

\noindent\Proof{} From Lemma~\ref{lem:L1}, the decomposition~\eqref{eq:decomp_L2} is stable by $\L.$ Moreover, the restriction of $\L$ to the space $\ker(P)$ is a sectorial operator whose spectrum is included in $\{z\in\Co \, |\, \Re e{(z)}>\eta\}$ for some small $\eta>0.$ The conclusion~\eqref{eq:lim_kerP} follows from~\cite{Hen81,Paz83}.\hfill$\Box$


\subsection{Proof of formula \eqref{eq:lim_tu_pushed}}\label{sec42}

Let $f$ satisfy the assumptions of Theorem~\ref{ref:theo_pushed} and let $(c,U)$ be either the pushed critical front when $c=c^*>2\sqrt{f'(0)}$ in case~$(A)$ or the unique front satisfying~\eqref{eq:pb_front} in cases~$(B)$ and~$(C).$ Let $\upsilon$ be the solution of the Cauchy problem~\eqref{eq:u0}-\eqref{eq:syst^i} and let $\tu$ be defined by $\tu(t,x)=\upsilon(t,x+ct).$ First of all, $\tu$ solves~\eqref{eq:pb_mv} and from the maximum principle the comparison~\eqref{ineq:tu_U} still holds. Moreover, since $U$ satisfies \eqref{eq:asymp_U_pushed}, \eqref{eq:asymp_U_bi} or \eqref{eq:asymp_U_ing}, $U$ and $\tu(t,\cdot)$ -- for all $t\geq0$ -- belong to the weighted space~$X_{c/2}$ defined by~\eqref{eq:def_X}. Next, let $\su$ be defined by~\eqref{eq:def_su*}. Since $\tu(t,\cdot)$ is in $X_{c/2}$ for all $t\geq0,$ the function $\su(t,\cdot)$ belongs to $L^2(\R),$ for all $t\geq0.$ Furthermore, $\varphi$ and $\su(0,\cdot)$ belong to $L^\infty(\R)$ from~\eqref{eq:asymp_U_pushed},~\eqref{eq:asymp_U_bi} and~\eqref{eq:asymp_U_ing}. From~\eqref{eq:pb_su*} and Lemma~\ref{lem:L1}, the initial condition $\su(0,\cdot)$ can be split in $L^2(\R)$ as follows:
\[
\su(0,\cdot)=P(\su(0,\cdot)) + w,
\]
where
\[\left\{\begin{array}{l}
\ds  P(\su(0,\cdot))=\lp\int_\R{\varphi(s)\,\su(0,s)\,ds}\rp \varphi \ \hbox{ and }\ \ \ds \varphi(x)=\frac{e^{cx/2}U(x)}{\ds \Big(\int_\R{U^2(s)\,e^{cs}\,ds}\Big)^{1/2}}\hbox{ for all }x\in\R, \\
 \ds  w=\su(0,\cdot)-P(\su(0,\cdot))\in\ker(P)\cap L^\infty(\R).
\end{array}\right.
\]
Since $\L\varphi=0,$ it follows that
\begin{equation}\label{eq:su_dec}
\su(t,\cdot)=P(\su(0,\cdot)) + e^{-t\L}w\hbox{ for all }t\ge0.
\end{equation}
Lemma~\ref{lem:L2} yields the existence of $C>0$ and $\eta>0$ such that
\begin{equation}\label{eq:su_P}
|\su(t,\cdot)-P(\su(0,\cdot))|_\infty=|e^{-t\L}w|_\infty\leq C e^ {-\eta t}|w|_\infty \hbox{ for all }t\ge0.
\end{equation}
Equation~\eqref{eq:su_dec} and the definition~\eqref{eq:def_su*} of $\su$ imply then that, for all $t>0$ and $x\in \R$,
\[\ds\begin{array}{rcl}
 \ds \tu(t,x) & = & \ds e^{-cx/2}\big( P(\su(0,\cdot))(x) + \big(e^{-t\L}w\big)(x)\big)\vspace{3pt}\\
           & = & \ds e^{-cx/2}\varphi(x) \lp\int_\R{\varphi(s)\,\su(0,s)\,ds}\rp + e^{-cx/2}\big(e^{-t\L}w\big)(x) \vspace{3pt}\\
           & = & p(\upsilon_0)U(x) + e^{-cx/2}\big(e^{-t\L}w\big)(x),
\end{array}
\]
where $p(\upsilon_0)\in(0,1]$ is given in~\eqref{eq:p}. It follows from~\eqref{eq:su_P} that $\upsilon(t,x+ct)-p(\upsilon_0)U(x)\to 0$ uniformly on compacts as $t\to+\infty$ and even uniformly in any interval of the type $[A,+\infty)$ with $A\in\R.$ This proves \eqref{eq:lim_tu_pushed}. \hfill$\Box$\break

Under an additional assumption on $\upsilon_0,$ the following lemma holds:

\begin{lem}\label{ref:lem_z}
Under the assumptions of Theorem~$\ref{ref:theo_pushed}$, if $\upsilon_0$  satisfies the additional assumption
\begin{equation}\label{eq:u0_pushed}
  \ds \limsup_{x\to-\infty}{\upsilon_0(x)}\leq p(\upsilon_0),\hbox{ or } \max{(\upsilon_0-p(\upsilon_0),0)}\in L^p(\R) \hbox{ for some } p\in[1,+\infty),
\end{equation}
then
\begin{equation}\label{eq:lim_sup_u_pushed}
  \sup_\R{\upsilon(t,\cdot)}\to p(\upsilon_0)  \hbox{ as } t\to+\infty.
\end{equation}
\end{lem}
\noindent{\it{Proof.}} The proof of~\eqref{eq:lim_sup_u_pushed}  is a consequence of~\eqref{eq:lim_tu_pushed} and Lemma~\ref{lem:sup_tu<lambda}. More precisely, let $\epsilon$ be any positive real number in $(0,1)$ and let $A$ be any real number. From the previous paragraph there is $t_0>0$ such that
\[
\upsilon(t,x)\leq p(\upsilon_0)U(x-ct) +\epsilon \leq p(\upsilon_0)+\epsilon \hbox{ for all }t\geq t_0\hbox{ and }x\ge A+ct.
\]
Lemma~\ref{lem:sup_tu<lambda} applied with $\lambda=\mu=p(\upsilon_0)$ implies that
$$\limsup_{t\to+\infty}{\Big(\sup_{x\in\R}{\upsilon(t,x)\Big)}}\leq p(\upsilon_0).$$
On the other hand, since $\liminf_{t\to+\infty}{\big(\sup_\R{\upsilon(t,\cdot)}\big)}\geq p(\upsilon_0)$ from~\eqref{eq:lim_tu_pushed} and $U(-\infty)=1,$ we get that $\sup_\R{\upsilon(t,\cdot)}\to p(\upsilon_0)$ as $t\to +\infty.$ \hfill$\Box$

\begin{remark}{\rm 
As in Remark~$\ref{rem:1}$, the above proof implies that if $\upsilon_0$ satisfies~\eqref{eq:u0} and $\upsilon_0(x)\to\lambda\in[0,1]$ as $x\to-\infty$, then $\sup_\R{\upsilon(t,\cdot)}\to\max{(\lambda,p(\upsilon_0))}$ as $t\to+\infty.$}
\end{remark}


\subsection{Spreading properties inside the pushed fronts: proofs of Theorem~\ref{ref:theo_pushed} and Proposition~\ref{ref:prop_pushed}}\label{sec:spread}

The previous section~\ref{sec42} shows that, in the pushed case, the right spreading speed of $\upsilon$ in the reference frame is equal to $c,$ in the sense that
$$c=\inf\big\{\gamma>0\ |\ \upsilon(t,\cdot+\gamma t)\to0\hbox{ uniformly in }(0,+\infty)\hbox{ as }t\to+\infty\big\}.$$
In this section, we prove that in the pushed case the left spreading speed of $\upsilon$ is actually at least equal to $0.$ More precisely, we prove that the solution $\upsilon$ moves to the left in the reference frame, at least at a sublinear rate proportional to $\sqrt{t},$ in the sense that $\liminf_{t\to+\infty}{\upsilon(t,\alpha\sqrt{t})}>0$ for all $\alpha\leq0$ (and, in fact, for all $\alpha\in\R$). This corresponds to formula \eqref{eq:sup_left_speed_pushed} in Theorem~\ref{ref:theo_pushed}. We also obtain some estimates, which are more precise than~\eqref{eq:lim_tu_pushed}, on the asymptotic profile of the solution $\upsilon$ in sets of the type $(\alpha\sqrt{t},+\infty)$ with $\alpha>0$ large enough. Lastly, we prove Proposition~\ref{ref:prop_pushed}, which shows that the solution $\upsilon$ cannot spread to the left with a positive speed if $\upsilon_0$ is  small near $-\infty,$ in the sense of~\eqref{eq:u0_pulled_lim}. The proofs of the pointwise estimates stated in Theorem~\ref{ref:theo_pushed} are based on formula~\eqref{eq:lim_tu_pushed} and on the construction of explicit sub- and super-solutions of~\eqref{eq:syst^i} in the reference frame.


\subsubsection{Description of the right spreading speed  and the asymptotic profile of solutions: proof of Theorem~\ref{ref:theo_pushed}}\label{sec431}

Let $f$ fulfill the assumptions of Theorem~\ref{ref:theo_pushed} and let $(c,U)$ be either the pushed critical front in the monostable case~$(A)$ with speed $c=c^*>2\sqrt{f'(0)}$, or the unique front in the bistable~$(B)$ and ignition~$(C)$ cases. First, as in the first inequality of~\eqref{ineq:j_A}, since $g(U(y))\to0^+$ as $y\to-\infty,$ there is $A\in\R$ such that
\begin{equation}\label{eq:g(U)A}
g(U(y))\geq 0\hbox{ for all } y\in(-\infty,A].
\end{equation}
Let $\upsilon$ be the solution of problem~\eqref{eq:syst^i} with initial condition $\upsilon_0$ satisfying~\eqref{eq:u0}. Theorem~\ref{ref:theo_pushed} implies that $\upsilon(t,A+ct)\to p(\upsilon_0)U(A)>0$ as $t\to+\infty.$ Choose any real number $\nu$ such that $0<\nu<p(\upsilon_0)U(A)$ and let $t_0>0$ be such that
\[
\upsilon(t,A+ct)\geq \nu \hbox{  for all } t\geq t_0.
\]

Now, let us construct a subsolution~$\underline{\upsilon}$ of the problem~\eqref{eq:syst^i} in the domain $t\geq t_0$ and $x\leq A+ct.$ Let $\underline{\upsilon}$ be defined by
$$\underline{\upsilon}(t,x)=h(t-t_0,x-A-ct_0)\hbox{ for all }t\ge t_0\hbox{ and }x\le A+ct,$$
where $h$ is the solution of the heat equation $\partial_t h=\partial^2_x h$ with initial condition $h(0,\cdot)=\nu \mathds{1}_{(0,+\infty)}$ in $\R$. On the one hand, there holds
$$\underline{\upsilon}(t,A+ct)\le\sup_{\R}{h(t-t_0,\cdot)}\leq \nu \leq \upsilon(t,A+ct) \hbox{ for all }t\geq t_0$$
and
$$ \underline{\upsilon}(t_0,x)= \nu \mathds{1}_{(0,+\infty)}(x-A-ct_0)=0\leq \upsilon(t_0,x)\hbox{ for all } x\in(-\infty,A+ct_0),$$
while, on the other hand, it follows from~\eqref{eq:g(U)A} that
$$\partial_t \underline{\upsilon}(t,x) - \partial^2_{x}\underline{\upsilon}(t,x)  -g(U(x-ct))\underline{\upsilon}(t,x)\leq 0$$
for all $t>t_0$ and $x\in(-\infty,A+ct)$. Then the maximum principle applied to~\eqref{eq:syst^i} implies that
\[
\upsilon(t,x)\geq \underline{\upsilon}(t,x)\hbox{ for all } t\geq t_0 \hbox{ and } x \le A+ct.
\]

Lastly, let $\alpha$ be any fixed real number. There exists $t_1>t_0$ such that, for all $t\ge t_1,$ one has $\alpha\sqrt{t}<ct+A$ and
\[
\upsilon(t,x)\geq \underline{\upsilon}(t,x)=h(t-t_0,x-A-ct_0)\hbox{ for all } x\in[\alpha\sqrt{t},ct+A].
\]
Since $h(t,\cdot)$ is increasing in $\R$ for all $t>0,$ we get, that for all $t\ge t_1$ and $x\in[\alpha\sqrt{t},ct+A],$
\[
\upsilon(t,x) \geq  \ds h(t-t_0,\alpha\sqrt{t} -A-ct_0) = \ds \frac{\nu}{\sqrt{\pi}}\ds \int_{(-\alpha\sqrt{t}+A+ct_0)/\sqrt{4(t-t_0)}}^{+\infty}{\ds e^{-y^2}\,dy}.
\]
Therefore,
\begin{equation*}
\liminf_{t\to+\infty}{\Big(\min_{\alpha\sqrt{t}\leq x\leq ct+A}{\upsilon(t,x)}\Big)}
          \geq \frac{\nu}{\sqrt{\pi}}\ds \int_{-\alpha/2}^{+\infty}{\ds e^{-y^2}dy} >0,
\end{equation*}
which together with~\eqref{eq:lim_tu_pushed} yields~\eqref{eq:sup_left_speed_pushed}.

Let us now turn to the proof of property~\eqref{eq:sup_left_speed_pushed2}. Let $\epsilon$ be any positive real number less than $2p(\upsilon_0)$. From~\eqref{eq:lim_tu_pushed} and $U\le U(-\infty)=1,$ Lemma~\ref{lem:bounds_u} applied with $\mu=p(\upsilon_0)$ (upper bound) and with $\mu=p(\upsilon_0)-\epsilon/2$ (lower bound) yields the existence of $\alpha_0>0,$ $A\in\R$ and $t_1>0$ such that $\alpha\sqrt{t}<A+ct$ for all $t\ge t_1$ and
\[
p(\upsilon_0)-\epsilon\leq \upsilon(t,x)\leq p(\upsilon_0)+\epsilon\ \hbox{ for all } t\geq t_1, \ x\in[\alpha\sqrt{t},A+ct] \hbox{ and } \alpha\geq \alpha_0.
\]
Even if it means decreasing $A,$ one can assume without loss of generality that $p(\upsilon_0)(1-U(A))\leq \epsilon.$ This implies that, for all $t\ge t_1$ and $\alpha\ge\alpha_0$,
\[
\ds \max_{\alpha\sqrt{t} \leq x \leq A+ct}{|\upsilon(t,x)-p(\upsilon_0)U(x-ct)|} \leq  \ds \max_{\alpha\sqrt{t}\leq x \leq A+ct}{|\upsilon(t,x)-p(\upsilon_0)|} + p(\upsilon_0)\big(1-U(A)\big) \leq  2\epsilon.
\]
From property~\eqref{eq:lim_tu_pushed} (and the fact that $\upsilon(t,x+ct)-p(\upsilon_0)U(x)\to 0$ as $t\to+\infty$ uniformly in $[A,+\infty)$, as observed at the end of the proof of~\eqref{eq:lim_tu_pushed}), we also know that there exists $t_2\geq t_1$ such that for all $t\geq t_2,$
\[
\max_{ x \geq A+ct}{|\upsilon(t,x)-p(\upsilon_0)U(x-ct)|}  \leq \epsilon.
\]
Remember that, in the above formula, the supremum is a maximum, from~\eqref{ineq:tu_U} and the continuity of $U$ and $\upsilon(t,\cdot)$ for all $t>0$. We conclude that
\[
\limsup_{t\to+\infty}{\Big(\max_{x\geq\alpha\sqrt{t} }{|\upsilon(t,x)-p(\upsilon_0)U(x-ct)|}\Big) } \leq 2\epsilon\,\hbox{ for all } \alpha\geq\alpha_0.
\]
Since $\epsilon>0$ is arbitrarily small, this completes the proof of Theorem~\ref{ref:theo_pushed}.\hfill$\Box$


\subsubsection{Description of the left spreading speed: proof of Proposition~\ref{ref:prop_pushed}}

In addition to~\eqref{eq:u0}, we assume that $\upsilon_0$ satisfies~\eqref{eq:u0_pulled_lim}. As already observed in the proof of Lemma~\ref{lem:sup_tu<lambda}, this implies that $\upsilon(t,x)\to0$ as $x\to-\infty$ for all $t>0.$ Let $\epsilon$ be any real number in $(0,1)$ and let $j_\epsilon$ be the function defined on $(-\infty,0)$ by~\eqref{eq:j_epsilon}. Let $A<0$ be such that~\eqref{ineq:j_A} holds. Since $\upsilon(1,x)\to0$ as $x\to-\infty,$ one can assume without loss of generality that $\upsilon(1,x)\leq\epsilon$ for all $x\leq A.$

As in the proof of Lemmas~\ref{lem:bounds_u} and~\ref{lem:sup_tu<lambda}, let us now construct a super-solution $\overline{\upsilon}$ of~\eqref{eq:syst^i} in the domain $t\geq1$ and $x\leq A.$ More precisely, let us set
\[
\overline{\upsilon}(t,x)=h(t-1,x-A)+j_\epsilon(x-ct) \hbox{ for all } t\ge 1 \hbox{ and } x\leq A,
\]
where $h$ solves the heat equation
\begin{equation*}
\left\{\begin{array}{rcll}
    \ds \partial_t h(t,x) & = & \partial_x^2 h(t,x), & t>0,\ x\in\R, \vspace{3pt}\\
    \ds h(0,x) & = & \epsilon\,\mathds{1}_{(-\infty,0)}(x)+2\,\mathds{1}_{(0,+\infty)}(x), & x\in\R.
\end{array}\right.
\end{equation*}
There holds $\upsilon(1,x)\le\epsilon\le\overline{\upsilon}(1,x)$ for all $x<A$, while $\upsilon(t,A)\le 1\le(2+\epsilon)/2=h(t-1,0)\le\overline{\upsilon}(t,A)$ for all $t>1$ from~\eqref{ineq:tu_U} and~\eqref{eq:j_epsilon}. Furthermore,
\begin{equation*}
\ds \partial_t \overline{\upsilon}(t,x)  -  \ds \partial_x^2 \overline{\upsilon}(t,x)  - g\big(U(x-ct)\big)\overline{\upsilon}(t,x)\geq 0 \hbox{ for all }  t>1 \hbox{ and } x<A,
\end{equation*}
as in~\eqref{eq:u_ou}. It follows from the maximum principle applied to~\eqref{eq:syst^i} that
$$\upsilon(t,x)\leq \overline{\upsilon}(t,x)\hbox{ for all }t\geq1\hbox{ and }x\leq A.$$
For any fixed $\alpha<0$, and $t>0$ large enough so that $\alpha\sqrt{t}<A,$ the maximum of $\upsilon(t,\cdot)$ on $(-\infty,\alpha\sqrt{t}]$ is reached since $\upsilon\geq0$ and $\upsilon(t,-\infty)=0$, and there holds
\[
\ds \max_{x\leq\alpha\sqrt{t}}{\upsilon(t,x)} \leq \ds \sup_{x\leq\alpha\sqrt{t}}{\big(h(t-1,x-A)+j_\epsilon(x-ct)\big)} \leq h(t-1,\alpha\sqrt{t}-A)+\epsilon
\]
since $h(t,\cdot)$ is increasing in $\R$ for all $t>0.$ Thus, for any fixed $\alpha<0$,
\[
\max_{x\leq\alpha\sqrt{t}}{\upsilon(t,x)} \leq 2\epsilon + \frac{2-\epsilon}{\sqrt{\pi}}\ds \int_{-(A+\alpha\sqrt{t})/\sqrt{4(t-1)}}^{+\infty}{\ds e^{-y^2}\,dy}
\]
for $t$ large enough, whence
\[
\limsup_{t\to+\infty}{\Big(\max_{x\leq\alpha\sqrt{t}}{\upsilon(t,x)}\Big)}\leq 3\epsilon \ \hbox{ for all }\alpha\le\alpha_0<0
\]
with $|\alpha_0|$ large enough. Since $\epsilon$ can be arbitrarily small, the proof of Proposition~\ref{ref:prop_pushed} is thereby complete.~\hfill$\Box$

\begin{remark}{\rm The proof of the lower bound of $\upsilon$ given in Subsection~\ref{sec431} implies immediately that, under the assumptions of Theorem~\ref{ref:theo_pushed}, $\liminf_{t\to+\infty}\upsilon(t,x)\ge p(\upsilon_0)/2$ locally uniformly in $x\in\R$. Furthermore, for any $\epsilon>0$, an adaptation of the above proof given in the present subsection implies that, under the assumptions of Proposition~\ref{ref:prop_pushed}, there are $A<0$ negative enough, $t_0>0$ positive enough and $B<\min(A+ct_0,0)$ negative enough such that $\upsilon(t,x)\le h(t-t_0,x)+j_{\epsilon}(x-ct)$ for all $t>t_0$ and $x\le ct+A$, where $h$ solves the heat equation $\partial_th=\partial^2_xh$ with initial condition $h(0,\cdot)=\epsilon\mathds{1}_{(-\infty,B)}+\mathds{1}_{(B,ct_0+A)}+(p(\upsilon_0)+\epsilon)\mathds{1}_{(ct_0+A,+\infty)}$. Therefore, since $\epsilon>0$ can be arbitrarily small, it follows that, under the assumptions of Proposition~\ref{ref:prop_pushed}, $\limsup_{t\to+\infty}\upsilon(t,x)\le p(\upsilon_0)/2$ locally uniformly in $x\in\R$, and finally $\upsilon(t,x)\to p(\upsilon_0)/2$ as $t\to+\infty$ locally uniformly in $x\in\R$.}
\end{remark}



\end{document}